\documentclass[a4paper,11pt]{article}

\usepackage{amssymb} 
\usepackage{amsmath}     
\usepackage{amsthm}
\usepackage{todonotes}
\usepackage{mathtools}
\usepackage{geometry}
\usepackage{a4wide}
\usepackage{graphicx}
\usepackage{hyperref}
\usepackage{subcaption}
\usepackage{float}
\usepackage[T1]{fontenc}

\usepackage[ruled,vlined,linesnumbered]{algorithm2e}
\SetKw{KwNot}{not}
\SetKw{KwOR}{or}
\SetKw{KwAND}{and}
\SetKw{KwExitFor}{exit for-loop}
\SetKw{KwStep}{Step}

\usepackage{authblk}

\newcommand{\Rset}{\mathbb{R}}

\newcommand{\OWA}{{\textup{OWA}}}
\newcommand{\OWAR}{{\textup{OWAR}}}
\newcommand{\new}[1]{\textcolor{black}{#1}}
\newcommand{\X}{{\mathcal{X}}}

\newtheorem{thm}{Theorem}

\newtheorem{cor}{Corollary}
\newtheorem{prop}{Proposition}

\newtheorem{definition}{Definition}

\begin{document}

\title{Robust Min-Max (Regret) Optimization using Ordered Weighted Averaging}

\author[1]{Werner Baak}
\author[1]{Marc Goerigk}

\affil[1]{Business Decisions and Data Science, University of Passau, Germany,

 \texttt{\{werner.baak,marc.goerigk\}@uni-passau.de}}

\author[2]{Adam Kasperski}

\affil[2]{Department of Operations Research
and Business Intelligence, Wroc{\l}aw University of Science and Technology, Poland,

\texttt{adam.kasperski@pwr.edu.pl}}

\author[3]{Pawe{\l} Zieli\'nski}

\affil[3]{Department of Fundamentals of Computer Science, Wroc{\l}aw University of Science and Technology, Poland,

\texttt{pawel.zielinski@pwr.edu.pl}}

\date{}

\maketitle

\begin{abstract}
In decision-making under uncertainty, several criteria have been studied to aggregate the performance of a solution over multiple possible scenarios. This paper introduces a novel variant of ordered weighted averaging (OWA) for optimization problems. It generalizes the classic OWA approach, which includes robust min-max optimization as a special case, as well as min-max regret optimization. We derive new complexity results for this setting, including insights into the inapproximability and approximability of this problem. In particular, we provide stronger positive approximation results that asymptotically improve the previously best-known bounds for the classic OWA approach. In computational experiments, we evaluate the quality of the proposed methods and compare the proposed setting with classic OWA and min-max regret approaches.
\end{abstract}



\section{Introduction}

In many real-world applications, decision-makers face uncertain and unpredictable scenarios that require careful consideration to reach optimal decisions. Uncertainty arises from various sources, such as incomplete or unreliable information, unforeseeable events, or unpredictable system dynamics. In this context, finding a good decision-making approach is crucial to address the consequences of uncertainty.
There are different methods available to tackle such problems, including \emph{stochastic optimization}
(see, e.g.,~\cite{birge2011introduction}) or \emph{robust optimization} (see  
surveys for~\cite{aissi2009min,goerigk2016algorithm} or for a guide~\cite{gorissen2015practical}).  In this paper, we focus on optimization problems with uncertainty in the objective function. This uncertainty is modeled by specifying a scenario set containing a finite number of cost realizations, called scenarios. In this context, the \emph{Ordered Weighted Averaging} (OWA for short)  criterion~\cite{YA88,YKB11} is commonly used.

 For the OWA approach, the idea to evaluate a solution is to sort its objective values over all possible scenarios, and to apply a weight vector to this sorted vector of values. The weights offer great flexibility to model preferences or risk-aversion of decision makers (see, e.g., \cite{baak2023preference,reimann2017well,xu2005overview}). 
 It turns out that many criteria used in decision-making under uncertainty, such as the maximum, average, median,  or Hurwicz (see,~e.g.,~\cite{LR57}), are special cases of OWA.  If we treat scenarios as a sample of random cost vectors, then OWA can be used to express the \emph{Conditional Value at Risk}~\cite{RU00} of this sample. The OWA criterion has also been used to aggregate objectives in a multi-objective optimization setting~\cite{OS03} or in problems where a feasible solution induces a multi-dimensional cost vector~\cite{CS19}.
 
 In general, the problem of minimizing OWA can be solved with the help of mixed-integer programming formulations~\cite{fernandez2014ordered,GS12,OO12,OS03}. The general case with arbitrary weights requires binary variables to express the ordering of the costs. However, the case of non-increasing weights is easier to handle, and more efficient models have become available \cite{CG15}. In particular, minimizing OWA is a convex problem if the underlying optimization problem is convex (for example, it is a linear programming one). The OWA criterion has also been applied to combinatorial problems, and some general results in this area have been shown in~\cite{KZ15}. Unfortunately, for most basic combinatorial problems (for example, for the shortest path, minimum spanning tree, or minimum assignment), minimizing OWA is NP-hard, even for two scenarios.
 Furthermore, for the general structure of weights, the problem is also not approximable. When the weights are non-increasing, a $O(K)$-approximation algorithm is known, where $K$ is the number of scenarios, provided that the underlying deterministic problem is polynomially solvable~\cite{KZ15}. This is the best general approximation algorithm known to date. An alternative approximation based on scenario aggregation and solving a MIP formulation has been proposed in~\cite{CGKZ20}.
 
In this paper, we generalize the classic OWA approach. We assume that each of the $K$ cost scenarios induces an affine function of a given solution. We then evaluate this solution by aggregating these affine functions using OWA.
 This approach also allows us to take into account the regrets of solutions under different scenarios. The \emph{min-max regret} approach has a long tradition in robust optimization~(see, e.g.,~\cite{KY97}). The maximum regret (also known as Savage~\cite{french1986decision,SA72}) criterion involves calculating the best possible outcome for each scenario and then considering the difference between the best and the actual outcome. In this paper, we provide a complexity characterization for the class of linear programming problems. We show that the problem of minimizing OWA is polynomially solvable when the weights are non-increasing or the number of scenarios is constant.
 On the other hand, we prove that the problem is strongly NP-hard and not approximable when the weights are non-decreasing. We also provide new results for the class of combinatorial problems.  We first establish some relationships between OWA minimization and $p$-norm minimization. In a $p$-norm minimization problem, a solution induces a $K$-dimensional vector of non-negative reals, and a $p$-norm is used to aggregate this vector into a single value~\cite{BCFM17}. By solving the $p$-norm minimization problem, we can strengthen the approximation results known to date.
 Using known results for $p$-norm minimization obtained in~\cite{BCFM17}, we apply this setting to derive new approximability results for some basic matroidal and network problems and provide a characterization of their approximability for various distributions of weights.

The remainder of this paper is structured as follows. In Section~\ref{sec:prel}, we recall the definition of the OWA criterion. We also show some known and new properties of OWA that are used later. In Section~\ref{sec:definition}, we formally introduce a general OWA optimization problem that we study and present some observations on its tractability. In Section~\ref{srlpowa}, we show that the case for non-decreasing weights is NP-hard, even if the nominal problem is a linear program. This proof also applies to the classic OWA setting. Our main results are presented in Section~\ref{sec:combinatorial}, where we discuss combinatorial problems. We provide inapproximability results in Section~\ref{sec:hardness} and new approximation results based on norm estimates in Section~\ref{sec:approx}. In Section~\ref{sec:appl}, we apply these results to the matroidal, the shortest path, and the minimum Steiner tree problems. 
We also discuss approximation guarantees stemming from scenario aggregation in Section~\ref{sec:aggregate}.
Moreover, in Section~\ref{sec:experiments}, we present computational experiments that compare our setting with classic min-max regret and OWA approaches, as well as experiments that evaluate the quality of our algorithms. We show that not only does the combination of OWA and regret provide a useful trade-off between the respective criteria, but it is also possible to find solutions that are close to optimality using our approximation algorithms.
In Section~\ref{sec:conclusion}, we conclude our paper.

\section{Preliminaries}
\label{sec:prel}

In this section, we recall the definition of the Ordered Weighted Averaging criterion. 
We also show some known and new inequalities, which are used later. We denote the set of nonnegative reals by $\Rset_{+}$ and $\overline{\Rset}_{+}=\Rset_{\geq 1}\cup\{\infty\}$ is the set of reals that are not smaller than~1 extended with $\infty$. We also use the notation $[K]=\{1,\dots,K\}$.
Let $\pmb{v}=(v_1,\dots,v_K)\in \Rset^K_{+} $ and $p\geq 1$. The value of
\[
\lVert \pmb{v}\rVert_p=\left( \sum_{i\in [K]} v^p_i\right)^{\frac{1}{p}}
\]
is called \emph{$p$-norm}. For $p=\infty$, we define $\lVert \pmb{v} \rVert_{\infty} = \max_{i\in[K]} v_i$ and  $1/p=0$.
We use the following well-known inequalities (see, e.g.,~\cite{HLP52,M70}):

\begin{prop}[H\"older's inequality]
\label{holin}
For every $\pmb{u},\pmb{v}\in \Rset^K_{+}$ and  $p,q\in \overline{\mathbb{R}}_{+}$ such that $\frac{1}{p}+\frac{1}{q}=1$,
the inequality
\[
  \sum_{i\in [K]}  u_i v_i\leq \lVert \pmb{u}\rVert_p \lVert \pmb{v}\rVert_q
\]
holds.
\end{prop}

\begin{prop}[Chebyshev's sum inequality]
 \label{czebin}
For every $\pmb{u},\pmb{v}\in \Rset^K$ such that $u_1\geq\cdots \geq u_K$ and
 $v_1\geq\cdots \geq v_K$  the inequality
 \[
 \left( \sum_{i\in [K]} u_i\right) \left( \sum_{i\in [K]} v_i\right) \leq K  \sum_{i\in [K]} u_iv_i
 \]
 holds.
\end{prop}

\begin{prop}[Rearrangement inequality]
\label{rearang}
	For every $\pmb{u},\pmb{v}\in \Rset^K_{+}$ such that $u_1\geq\cdots \geq u_K$ and
 $v_1\geq\cdots \geq v_K$ and any permutation $\pi$ of $[K]$  the following inequality
 $$\sum_{i\in [K]} u_i v_i\geq \sum_{i\in [K]} u_i v_{\pi(i)}$$
 holds.
\end{prop}

\begin{prop}
\label{norinq}
For every $\pmb{v}\in \Rset^K_{+}$ and $p,q\in\overline{\Rset}_{+}$ such that $p\leq q$,
the inequalities
\[
  \lVert \pmb{v}\rVert_q \leq \lVert \pmb{v}\rVert_p\leq
  \frac{K^{\frac{1}{p}}}{K^{\frac{1}{q}}}
  \lVert\pmb{v}\rVert_q
\]
hold.
\end{prop}

Let us recall the definition of the Ordered Weighted Averaging (OWA for short).

\begin{definition}[\cite{YA88}]
Let $\pmb{w}=(w_1,\dots,w_K)\in \Rset_{+}^K$ be a vector of nonnegative weights such that  $w_1+w_2+\dots+w_K=1$. The Ordered Weighted Averaging of $\pmb{v}\in \Rset^K$ is defined as
$$\OWA_{\pmb{w}}(\pmb{v})=\sum_{i\in [K]} w_i v_{\pi(i)},$$
where $\pi$ is a permutation of~$[K]$ such that
$v_{\pi(1)}\geq v_{\pi(2)}\geq \dots \geq v_{\pi(K)}$.
\end{definition}
Let us describe several special cases of OWA.
If $w_1=1$ and $w_i=0$ for $i=2,\dots,K$, then OWA becomes the maximum.
 If $w_K=1$ and $w_i=0$ for $i=1,\dots,K-1$, then OWA becomes the minimum. In general, if $w_k=1$ and $w_i=0$ for $i\in [K]\setminus\{k\}$, then OWA is the $k$-th largest
component of~$\pmb{v}$.
 In particular, when $k=\lfloor K/2 \rfloor +1$, the $k$-th largest component is the median.
If $w_i=\frac{1}{K}$ for all $i \in [K]$, i.e. when the weights are \emph{uniform}, then OWA is the average. Finally, if $w_1=\lambda$ and $w_K=1-\lambda$, for some fixed $\lambda\in [0,1]$, and $w_i=0$ for the remaining weights, then we get the convex combination of the maximum and the minimum components of $\pmb{v}$ (in decision-making it is called the Hurwicz criterion).

In this paper, we mainly discuss the case of \emph{non-increasing weights}, i.e. when $w_1\geq w_2\geq \dots\geq w_K$. For this particular structure of weights, OWA is also called an \emph{ordered norm}~\cite{CS19}. The case of non-increasing weights can lead to problems easier from the computational point of view due to the following fact:
\begin{prop}
\label{propconv}
If $\pmb{w}$ is non-increasing, then $\OWA_{\pmb{w}}(\pmb{v})$ is  a convex function in $\Rset^{K}_+$.
\end{prop}
\begin{proof}
For each $\lambda\in [0,1]$ and any two vectors $\pmb{u},\pmb{v}\in \Rset^{K}_+$
$$\OWA_{\pmb{w}}(\lambda\pmb{u}+(1-\lambda)\pmb{v})=\sum_{i\in [K]} w_i (\lambda u_{\pi(i)}+(1-\lambda)v_{\pi(i)})=\lambda\sum_{i\in [K]} w_i u_{\pi(i)}+(1-\lambda)\sum_{i\in [K]} w_i v_{\pi(i)},$$
where $\lambda u_{\pi(1)}+(1-\lambda)v_{\pi(1)}\geq \dots \geq \lambda u_{\pi(K)}+(1-\lambda)v_{\pi(K)}$.
Proposition~\ref{rearang} yields $\OWA_{\pmb{w}}(\pmb{u})\geq \sum_{i\in [K]} w_i u_{\pi(i)}$ and $\OWA_{\pmb{w}}(\pmb{v})\geq \sum_{i\in [K]} w_i v_{\pi(i)}$ which implies
$$\OWA_{\pmb{w}}(\lambda\pmb{u}+(1-\lambda)\pmb{v})\leq \lambda \OWA_{\pmb{w}}(\pmb{u}) + (1-\lambda) \OWA_{\pmb{w}}(\pmb{v}).$$
\end{proof}

We now prove the following estimates on the OWA value.

\begin{prop}
For every vector $\pmb{v}\in \Rset^K_{+}$ and non-increasing weights $\pmb{w}$, the inequalities
\begin{equation}
\label{el3e}
\mathrm{OWA}_{\pmb{w}}(\pmb{v})\leq  \lVert \pmb{w}\rVert_{q}  \lVert \pmb{v}\rVert_p
\leq \rho
\, \mathrm{OWA}_{\pmb{w}}(\pmb{v})
\end{equation}
hold,
where $p, q\in\overline{\mathbb{R}}_{+} $, $\frac{1}{p}+\frac{1}{q}=1$,
and $\rho= K^{\frac{1}{p}}\frac{ \lVert \pmb{w}\rVert_{q} }{\lVert \pmb{w}\rVert_p}$.
\label{eglp}
\end{prop}
\begin{proof}
The first inequality in~(\ref{el3e}) follows from Proposition~\ref{holin}, namely
\begin{align*}
\mathrm{OWA}_{\pmb{w}}(\pmb{v})&=\sum_{i \in [K]} w_i v_{\pi(i)}\leq
\left( \sum_{i\in [K]} w^q_i\right)^{\frac{1}{q}}\left( \sum_{i\in [K]} v_{\pi(i)}^p\right)^{\frac{1}{p}}=
\left( \sum_{i\in [K]} w^q_i\right)^{\frac{1}{q}}\left( \sum_{i\in [K]} v_{i}^p\right)^{\frac{1}{p}}\\
&=\lVert \pmb{w}\rVert_q  \lVert \pmb{v}\rVert_p.
\end{align*}
The second inequality in~(\ref{el3e})
holds as
\begin{align*}
\lVert \pmb{w}\rVert_q  \lVert \pmb{v}\rVert_p&=
\frac{\lVert \pmb{w}\rVert_q }{\lVert \pmb{w}\rVert_p} \lVert\pmb{w}\rVert_p \lVert\pmb{v}\rVert_p
= \frac{\lVert \pmb{w}\rVert_q}{\lVert \pmb{w}\rVert_p}
\left( \sum_{i\in [K]} w^p_{i}\right)^{\frac{1}{p}}\left( \sum_{i\in [K]} v_{\pi(i)}^p\right)^{\frac{1}{p}}\\
&\leq K^{\frac{1}{p}} \frac{\lVert \pmb{w}\rVert_q}{\lVert \pmb{w}\rVert_p}
\left( \sum_{i\in [K]} w^p_{i}  v_{\pi(i)}^p\right)^{\frac{1}{p}}
\leq K^{\frac{1}{p}} \frac{\lVert \pmb{w}\rVert_q}{\lVert \pmb{w}\rVert_p}
\left( \sum_{i\in [K]} w_{i}  v_{\pi(i)}\right)=
K^{\frac{1}{p}} \frac{\lVert \pmb{w}\rVert_q}{\lVert \pmb{w}\rVert_p}
\, \mathrm{OWA}_{\pmb{w}}(\pmb{v}),
\end{align*}
where we used Proposition~\ref{czebin} for the first inequality and Proposition~\ref{norinq} for the second one.
\end{proof}

Finally, we recall the following result from the literature.

\begin{prop}[\cite{CGKZ20}]
\label{lemaggr}
Assume that $K$ is a multiple of $\ell\geq 1$. For each $\pmb{v}\in \Rset_{+}^K$ and non-increasing weights $\pmb{w}$ the inequalities
\begin{equation}
\OWA_{\overline{\pmb{w}}}(\overline{\pmb{v}}) \le \OWA_{\pmb{w}}(\pmb{v}) \le \ell\phi\OWA_{\overline{\pmb{w}}}(\overline{\pmb{v}})
\end{equation}
hold for $\phi = \max_{k\in[K/\ell]} \left(\sum_{i\in[k]} w_i / \sum_{i\in[k]} \overline{w}_i\right)$, where $\overline{\pmb{v}} = (\overline{v}_1,\ldots,\overline{v}_{K/\ell})\in \Rset_{+}^{K/\ell}$
with $\overline{v}_i = (v_{(i-1)\ell+1} + \ldots + v_{i\ell})/\ell$  and
$\overline{\pmb{w}} = (\overline{w}_1,\ldots,\overline{w}_{K/\ell}) \in \Rset_{+}^{K/\ell}$
with $\overline{w}_i = w_{(i-1)\ell+1} + \ldots + w_{i\ell}$ for $i\in[K/\ell]$.
\end{prop}

\section{Problem formulation}
\label{sec:definition}

Consider the following generic optimization problem~$\mathcal{P}$:
$$\mathcal{P}:\; \min_{\pmb{x}\in \mathcal{X}} \pmb{c}^T\pmb{x},$$
where $\pmb{x}=(x_1,\dots,x_n)^T$ is a vector of decision variables, $\mathcal{X}$ is a set of feasible solutions and $\pmb{c}\in \Rset^n_{+}$ is a given cost vector. The set $\mathcal{X}$ is typically described by a system of linear constraints involving the variables $x_1,\dots,x_n$. If $\mathcal{X}$ is a polyhedron in~$\Rset_{+}^n$, then $\mathcal{P}$ is a linear programming problem. If additionally $\mathcal{X}\subseteq \{0,1\}^n$, then $\mathcal{P}$ is a combinatorial optimization problem. In particular, we get an important class of network problems assuming that $\mathcal{X}$ is a set of characteristic vectors of some objects in a given graph. For example, $\mathcal{P}$ can be the shortest path, the minimum spanning tree,
the minimum assignment problem, etc.~\cite{AMO93}.

In many practical situations, the cost vector $\pmb{c}$ is uncertain,
which means that precise values of its components are not known before a solution is to be computed.
In this case, a \emph{scenario set}~$\mathcal{U}$ containing possible realizations of $\pmb{c}$ is a
part of the input. Each realization in~$\mathcal{U}$ is called \emph{scenario}, corresponding to a possible state of the world.
 In this paper, we use the \emph{discrete uncertainty representation}~\cite{KY97}, namely,
 the uncertainty set~$\mathcal{U}=\{\pmb{c}_1,\dots,\pmb{c}_K\}\subset \mathbb{R}^n_{+}$ contains a finite number of $K\geq 1$ scenarios.  The scenarios can be listed explicitly or can result from a sampling of the uncertain (random) cost vector.
In this paper we use the convenient representation of~$\mathcal{U}$ as a
 scenario matrix $\pmb{C}\in \Rset^{K\times n}_{+}$,   where the $i$-th row $\pmb{c}_i$ of $\pmb{C}$ is the $i$-th scenario $\pmb{c}^T_i=(c_{i1},\ldots,c_{in})\in \Rset^{n}_{+}$, $i\in[K]$.  Let $\pmb{b}\in \Rset^K_{+}$ be a fixed vector. Define
$$\pmb{v}_{\pmb{b}}(\pmb{x})=\pmb{C}\pmb{x}-\pmb{b}=(\pmb{c}^T_1\pmb{x}-b_1,\ldots,\pmb{c}^T_K\pmb{x}-b_K)^T.$$
We assume that $\pmb{c}_i^T\pmb{x}-b_i\geq 0$ for any solution $\pmb{x}\in \mathcal{X}$ and each $i\in [K]$, so $\pmb{v}_{\pmb{b}}(\pmb{x})\in \Rset_{+}^K$ for each feasible solution $\pmb{x}\in\mathcal{X}$.  In particular, $b_i$ can be the optimal objective value of~$\mathcal{P}$ under scenario $\pmb{c}_i$ or a lower bound on this value. Clearly,
if $\pmb{b}=\pmb{0}$, then $\pmb{v}_{\pmb{b}}(\pmb{x})$ is a vector of the solution costs under scenarios $\pmb{c}_1,\dots,\pmb{c}_K$. Using an additional vector $\pmb{b}$, whose components~$b_i$
are  the optimal objective values of~$\mathcal{P}$ under scenarios $\pmb{c}_i$, $i\in [K]$, respectively,
allows us to express the regrets of solution $\pmb{x}$, namely $\pmb{c}_i^T\pmb{x}-b_i$ can be interpreted as a regret of solution $\pmb{x}$ under scenario $\pmb{c}_i$.
Given a vector of weights $\pmb{w}$ we wish to investigate
the following optimization problem:
\begin{equation}
      \textsc{Owa}~\mathcal{P}:\;\min_{\pmb{x}\in\mathcal{X}}\mathrm{OWA}_{\pmb{w}}(\pmb{v}_{\pmb{b}}(\pmb{x})).
      \label{owap}
\end{equation}
The OWA criterion used in~(\ref{owap}) aggregates $K$ affine functions, in particular the costs or regrets of solution~$\pmb{x}$ under scenarios $\pmb{c}_1,\dots,\pmb{c}_K$.
Note that in the classic definition of OWA optimization (see, e.g., \cite{chassein2020approximating,fernandez2014ordered}), we have $\pmb{b}=\pmb{0}$.
Problem~(\ref{owap}) thus encompasses a broader family of optimization problems, which use various criteria for decision-making under uncertainty (see, e.g.,~\cite{LR57}).

\begin{figure}[ht]
	\centering
	\includegraphics[height=3.2cm]{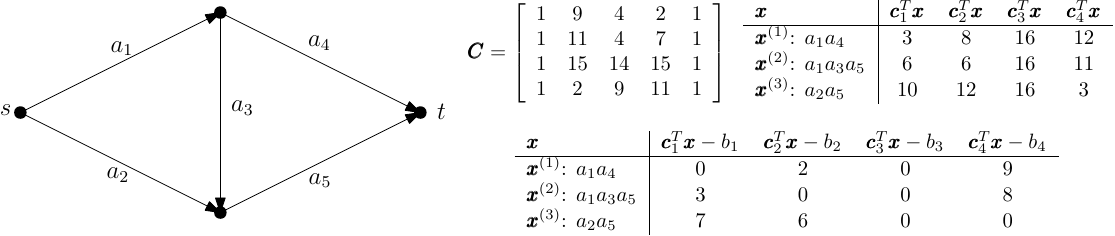}
	\caption{A sample shortest path problem with $K=4$ cost scenarios.} \label{figex}
\end{figure}

Let us illustrate the problem using a small example. Consider a network shown in Figure~\ref{figex} in which we seek a shortest path from node $s$ to node $t$. There are $K=4$ cost scenarios being the rows of matrix $\pmb{C}$. The problem has exactly three solutions: $\pmb{x}^{(1)}$, $\pmb{x}^{(2)}$, $\pmb{x}^{(3)}$ being the characteristic vectors of the paths $a_1a_4$, $a_1a_3a_5$ and $a_2a_5$, respectively. In Figure~\ref{figex} the costs $\pmb{c}_i^T\pmb{x}$ and the regrets $\pmb{c}_i^T\pmb{x}-b_i$ for each path are shown. These quantities are the components of the vector $\pmb{v}_{\pmb{b}}(\pmb{x})$ for $\pmb{b}=(0,0,0,0)^T$ and  $\pmb{b}=(3,6,16,3)^T$, respectively.

If $\pmb{b}=\pmb{0}$ and we pay attention only to the largest path cost, then all three solutions are equivalent, as the largest cost of each path is the same and occurs in scenario $\pmb{c}_3$. We can observe here the so-called \emph{drowning effect}~\cite{DF99} in which only one bad scenario is taken into account, and the information contained in other scenarios is ignored. This drawback worsens as the number of scenarios increases because the worst scenario can become less and less likely. Therefore, choosing some weight vector and using the OWA criterion for evaluating solutions is reasonable. For example, if $\pmb{w}=(0.6, 0.3, 0.1,0)$, then the solution $\pmb{x}^{(2)}$ is optimal with $\OWA_{\pmb{w}}(\pmb{v}_{\pmb{0}}(\pmb{x}^{(2)})) = 0.6\cdot 16 + 0.3\cdot 11 + 0.1\cdot 6+0\cdot 6 = 13.5$. Observe also that for $\pmb{b}=\pmb{0}$ and for any vector of non-increasing weights, the solution $\pmb{x}^{(1)}$ is not worse than $\pmb{x}^{(3)}$. If $\pmb{b}=(3,6,16,3)^T$ and we pay attention only to the largest path regret (the largest opportunity loss), then we should choose the solution $\pmb{x}^{(3)}$ whose maximum regret is equal to~7. Observe, however, that the second largest regret of $\pmb{x}^{(3)}$ is~6, and this information is ignored. Therefore, using the OWA criterion with some weight vector can also be reasonable. If we again use $\pmb{w}=(0.6, 0.3, 0.1,0)$, then the optimal solution is $\pmb{x}^{(2)}$ with 
$\OWA_{\pmb{w}}(\pmb{v}_{\pmb{b}}(\pmb{x}^{(2)})) = 0.6\cdot 8 + 0.3\cdot 3 + 0.1\cdot 0+0\cdot 0 = 5.7$.

In more detail, let us investigate some special cases of $\textsc{Owa}~\mathcal{P}$.
If $\pmb{w}=(1,0,\dots,0)$, then we get the following robust min-max (regret) problem~\cite{KY97}:
$$\textsc{Min-max}~\mathcal{P}: \min_{\pmb{x}\in \mathcal{X}} \lVert \pmb{v}_{\pmb{b}}(\pmb{x}) \rVert _{\infty}.$$
 On the other hand, the uniform weight vector $\pmb{w}=(\frac{1}{K},\dots,\frac{1}{K})$ leads to  the following problem:
 $$\textsc{Min-average}~\mathcal{P}: \; \min_{\pmb{x}\in \mathcal{X}}\frac{1}{K} \lVert \pmb{v}_{\pmb{b}}(\pmb{x}) \rVert _{1}. $$
 It is easy to check that  $\textsc{Min-average}~\mathcal{P}$ can be reduced to solving the deterministic problem~$\mathcal{P}$ with the cost vector $\hat{\pmb{c}}=\frac{1}{K}\sum_{i\in [K]} \pmb{c}_i$, which follows immediately from the fact that in this case
 $$\OWA_{\pmb{w}}(\pmb{v}_{\pmb{b}}(\pmb{x}))=\frac{1}{K}\sum_{i\in [K]} (\pmb{c}_i^T\pmb{x}-b_i)=\hat{\pmb{c}}^T\pmb{x}-\frac{1}{K}\sum_{i\in [K]} b_i$$
 and $b_i$, $i\in [K]$, are constant.

 It is worth pointing out that
  the $\textsc{Min-max}~\mathcal{P}$ and $\textsc{Min-Average}~\mathcal{P}$ problems are boundary cases of $\textsc{Owa}~\mathcal{P}$ with non-increasing weights. Generally, a non-increasing weight vector $\pmb{w}$ can be used to model risk-averse decision-makers. The less uniform $\pmb{w}$ is, the more risk-averse the decision maker is. In the boundary case $\pmb{w}=(1,0,\dots,0)$, the decision maker is extremely risk averse and pays attention only to the worst scenario that can occur for solution $\pmb{x}$.  On the other hand, if the weights are uniform, the decision maker is risk-neutral.
  We get another important special case of non-increasing weights by choosing $\pmb{w}=(\frac{1}{l},\dots,\frac{1}{l},0,\dots,0)$, for $l\in [K]$, where only the first $l$ weights in $\pmb{w}$ are positive. It is easy to see that $\OWA_{\pmb{w}}(\pmb{v}_{\pmb{b}}(\pmb{x}))$ is then the average of the $l$-largest values in $\pmb{v}_{\pmb{b}}(\pmb{x})$ and $\textsc{Owa}~\mathcal{P}$ is a special case of the problem of minimizing the \emph{Conditional Value at Risk} (CVaR for short) (see~\cite{RU00}). Indeed, if we interpret~$(\pmb{c}_i, b_i)_{i\in [K]}$ as a sample of some random vector $(\tilde{\pmb{c}},\tilde{b})$, then
 \begin{equation}
 \label{cvar}
 \OWA_{\pmb{w}}(\pmb{v}_{\pmb{b}}(\pmb{x}))=\inf\left\{t \in \Rset\,:\, t+\frac{1}{\alpha K}\sum_{i\in [K]} [\pmb{c}_i^T\pmb{x}-b_i-t]_{+}\right\},
 \end{equation}
 where $\alpha=\frac{1}{l}$ and $[u]_{+}=\max\{0,u\}$.
  The right-hand side of~(\ref{cvar}) is the Conditional Value at risk of the sample $(\pmb{c}_i, b_i)_{i\in [K]}$ with a risk level of $\alpha$ (see, e.g.,~\cite{P00, RU00}).

  From Proposition~\ref{propconv}, we have the following result:
  \begin{prop}
\label{convp}
	If $\pmb{w}$ is non-increasing and $\mathcal{X}$ is a convex set, then $\textsc{Owa}~\mathcal{P}$ is a convex optimization problem.
\end{prop}
Hence, $\textsc{Owa}~\mathcal{P}$ for non-increasing weights can be a tractable optimization
problem (see, e.g.,~\cite{BV04}).
If $\pmb{w}$ is non-increasing, then using Proposition~\ref{rearang}, we get
\begin{equation}
\label{owaperm}
 \OWA_{\pmb{w}}(\pmb{v}_{\pmb{b}}(\pmb{x})) = \max_{\pi\in\Pi} \sum_{i\in[K]} w_i (\pmb{c}^T_{\pi(i)}\pmb{x} -b_{\pi(i)}),
 \end{equation}
where $\Pi$ denotes the set of all permutations of~$[K]$. Representing $\Pi$ by the assignment constraints and using the dual of~(\ref{owaperm}) leads to the following compact reformulation of $\textsc{Owa}~\mathcal{P}$
(see~\cite{CG15}):
\begin{equation}
\label{mipowa1}
	\begin{array}{lll}
		\min & \displaystyle \sum_{k\in[K]} (\alpha_k + \beta_k) \\
			\text{s.t. } & \displaystyle \alpha_k + \beta_i \ge w_i (\pmb{c}_k^T\pmb{x}-b_k) & \forall i,k\in[K] \\
& \pmb{x}\in\mathcal{X}
	\end{array}
\end{equation}
Observe that~(\ref{mipowa1}) is a convex problem if $\mathcal{X}$ is a convex set.

  For arbitrary weight vector $\pmb{w}$, the $\textsc{Owa}~\mathcal{P}$ problem can be represented as the following program~\cite{OO12}:
\begin{equation}
\label{mipowa}
\begin{array}{llll}
	\min  & \displaystyle \sum_{i\in [K]} w_i y_i \\
		 &\displaystyle y_k+  M z_{ik} \geq \pmb{c}_i^T\pmb{x}-b_i & \forall i,k\in [K]\\
		 & \displaystyle\sum_{i\in [K]} z_{ik}\leq k-1 & \forall k\in [K] \\
		 & \displaystyle z_{ik}\in \{0,1\} & \forall  i,k\in [K]\\
		 & \pmb{x}\in \mathcal{X}
\end{array}
\end{equation}
where $M\ge \max_{\pmb{x}\in\mathcal{X}} \max_{i\in[K]} \pmb{c}_i^T \pmb{x}$ is a sufficiently large constant. Notice that model~(\ref{mipowa}) has $K^2$ binary variables that express the ordering of $\pmb{c}_i^T\pmb{x}-b_i$ for $i\in [K]$.

 We now show some methods of solving~$\textsc{Owa}~\mathcal{P}$, which can be used for a particular structure of weights, which do not need to be non-increasing. Assume first that $\pmb{w}=(\lambda,0,\dots,0,1-\lambda)$ for $\lambda\in [0,1]$, so OWA  is the Hurwicz criterion, being a convex combination of the maximum and the minimum value of the solution $\pmb{x}$. The $\textsc{Owa}~\mathcal{P}$ problem can be then rewritten as follows:
\begin{equation}
\label{hur}
\min_{\pmb{x}\in \mathcal{X},k\in [K]}\left(\lambda \max_{i\in [K]} (\pmb{c}_i^T\pmb{x}-b_i)+(1-\lambda)(\pmb{c}_k^T\pmb{x}-b_k)\right).
\end{equation}
An optimal solution to~(\ref{hur}) can be found by solving $K$ $\textsc{Min-max}~\mathcal{P}$ problems, i.e. the problems with $\pmb{w}'=(1,0,\dots,0)$. By Proposition~\ref{convp}, we get a tractable problem when the set $\mathcal{X}$ is convex. We now show a generalization of~(\ref{hur}). Suppose that the weight vector $\pmb{w}=(w_1,\dots,w_{r},0,\dots,0,w_{s},\dots,w_{K})$, where $w_1\geq\dots\geq w_r\ge 0$ and $0\le w_s\leq \dots \leq w_K$. Therefore, the first $r$ weights are non-increasing, while the last $K-s+1$ weights are non-decreasing. The OWA criterion with $\pmb{w}$ can be seen as a generalization of the Hurwicz criterion. Let $\pmb{w}'=(w_1,\dots,w_r,0,\dots,0)$. Notice that $\pmb{w}'$ is non-increasing. Let $\mathcal{S}$ be the set of all permutations of all $(K-s+1)$-element subsets of $K$. The $\textsc{Owa}~\mathcal{P}$ problem can be then expressed as
\begin{equation}
\label{hurgen}
\min_{\pmb{x}\in \mathcal{X},\sigma\in \mathcal{S}} \left(\OWA_{\pmb{w}'}(\pmb{v}_{\pmb{b}}(\pmb{x}))+\sum_{i=s}^K w_i(\pmb{c}_{\sigma(i-s+1)}^T\pmb{x}-b_{\sigma(i-s+1)})\right).
\end{equation}
Again, if $\mathcal{X}$ is convex, then an optimal solution to~(\ref{hurgen}) can be found by solving a family of $|\mathcal{S}|$ convex problems (notice that for fixed $\sigma$ the problem~(\ref{hurgen}) is convex as the objective function is a sum of convex functions). We get an efficient algorithm only when the size of $\mathcal{S}$ is not large, which means that only several of the last weights in $\pmb{w}$ are positive. This is, in particular, the case if $s$ is a constant number. In Section~\ref{srlpowa}, we show that $\textsc{Owa}~\mathcal{P}$ with non-decreasing weights can be intractable even if $\mathcal{X}$ is a convex set.

\section{Robust linear programming with the OWA criterion}
\label{srlpowa}

In this section, we assume that $\mathcal{X}$ is a polytope in $\Rset_{+}^n$, i.e. a closed and bounded subset of~$\Rset_{+}^n$, which can be described by a system of linear constraints on the real variables $x_1,\dots,x_n$. Because $\mathcal{X}$ is convex, Proposition~\ref{convp} implies that for non-increasing weights, the problem can be solved by~(\ref{mipowa1}), which is a linear programming problem. Furthermore, $\textsc{Owa}~\mathcal{P}$ can be solved in polynomial time for any vector of weights if $K$ is constant. Indeed, for a constant $K$, the formulation~(\ref{mipowa}) has $K^2$ binary variables, which is also constant. We can thus solve~(\ref{mipowa}) by trying all possible assignments to the binary variables. The following result characterizes the problem complexity when $K$ is part of the input and the vector of weights is non-decreasing.

\begin{thm}
	If $K$ is part of the input, then $\textsc{Owa}~\mathcal{P}$ is strongly NP-hard and not approximable 
	unless P=NP,
	if the vector of weights $\pmb{w}$ is non-decreasing, $\pmb{b}=\pmb{0}$, and $\mathcal{X}$ is a polytope in $[0,1]^n$.
\end{thm}
\begin{proof}
Consider the following \textsc{MINSAT} problem. Given a set of $s$  boolean variables $q_1,\dots q_s$, a collection of $t$ clauses $\mathcal{C}_1,\dots, \mathcal{C}_t$ over the boolean variables and a positive integer $r< t$. 
We ask if there is a truth assignment to the variables in which at least $r$ clauses are unsatisfied. The \textsc{MINSAT} problem is known to be strongly NP-complete, even if each clause contains at most two literals~\cite{KM94}.  

Given an instance of   \textsc{MINSAT}, we build the corresponding instance of $\textsc{Min-owa}~\mathcal{P}$ as follows. 
Let us define variables $x_i$ and $\overline{x}_i$ for each $i\in [s]$, so the number of variables is $n=2s$. 
Define the polytope $\mathcal{X}\subset [0,1]^n$ by  constraints $x_i\geq 0$, $\overline{x}_i\geq 0$ and $x_i+\overline{x}_i=1$, $i\in [s]$. For each clause $\mathcal{C}_i$, we form scenario $\pmb{c}_i$ as follows.
If $q_j\in \mathcal{C}_i$, then the cost of $x_j$ is~1; if $\overline{q}_j\in \mathcal{C}_i$
($\overline{q}_j$ is the negation of~$q_j$), then the cost of $\overline{x}_j$
is~1 under $\pmb{c}_i$; the costs of the remaining variables under $\pmb{c}_i$ are set to~0. The number of scenarios is $K=t$.  We set $\pmb{b}=\pmb{0}$. The vector of weights is

$$\pmb{w}=\left(\underbrace{0,\dots,0}_{K-r\; \text{weights}}, \underbrace{\frac{1}{r}, \dots, \frac{1}{r}}_{r \;\text{positive weights}} \right)$$

 To illustrate the reduction, consider a sample instance with variables $q_1,q_2,q_3,q_4$, clauses $(q_1 \vee \overline{q}_2)$,
  $(\overline{q}_2\vee q_3)$, $(q_1 \vee q_3)$, $(\overline{q}_3 \vee q_4)$, 
  $(q_1\vee \overline{q}_4)$, $(\overline{q}_2 \vee \overline{q}_4)$, $(\overline{q}_1,\vee q_2)$, $(q_1 \vee q_2)$ and $r=3$.  
  The scenarios for this instance are shown in Table~\ref{tab1} and $\pmb{w}=(0,0,0,0,0,\frac{1}{3},\frac{1}{3}, \frac{1}{3})$.
 \begin{center}
 \begin{table}[ht]
 	 \centering
 	 \caption{Scenarios for the sample instance of \textsc{MINSAT}.}\label{tab1}
 	\begin{tabular}{cc|cccccccccccccccc}
    	  $\pmb{x}$ & $\pmb{c}_1$ & $\pmb{c}_2$ & $\pmb{c}_3$ & $\pmb{c}_4$ & $\pmb{c}_5$ & $\pmb{c}_6$ & $\pmb{c}_7$ & $\pmb{c}_8$ \\ \hline
	 $x_1$ & 1 & 0 & 1 & 0 & 1 & 0 & 0 & 1 \\
	 $\overline{x}_1$ & 0 & 0 & 0 & 0 & 0 & 0 & 1 & 0 \\
	 $x_2$ & 0 & 0 & 0 & 0 & 0 & 0 & 1 & 1  \\
	 $\overline{x}_2$ & 1 & 1 & 0 & 0 & 0 &1 & 0 & 0  \\
	 $x_3$ & 0 & 1 & 1 & 0 & 0 & 0 & 0 & 0 \\
	 $\overline{x}_3$ & 0 & 0 & 0 & 1 & 0 & 0 & 0 & 0 \\
	 $x_4$ & 0 & 0 & 0 & 1 & 0  & 0 & 0 & 0 \\
	 $\overline{x}_4$ & 0 & 0 & 0 & 0 & 1 & 1 & 0 & 0\\
	\end{tabular}
 \end{table}
 \end{center}

 We now show that the answer to \textsc{MINSAT} is \emph{yes} if  and only if $\OWA_{\pmb{w}}(\pmb{v}_{\pmb{b}}(\pmb{x}))\leq 0$ for some feasible solution $\pmb{x}\in\mathbb{X}$.
  
 Assume that the answer to \textsc{MINSAT} is \emph{yes} and let $q_1,\dots,q_s$ be a truth assignment to the variables for which at least $r$ clauses are not satisfied. Let us form a feasible solution $\pmb{x}\in \mathcal{X}$ such that $x_j=1, \overline{x}_j=0$ if $q_j=1$  and $x_j=0, \overline{x}_j=1$ if $q_j=0$. By the construction, there are at least $r$ scenarios under which the cost of $\pmb{x}$ is~0 and thus no more than $K-r$ scenarios under which the cost of $\pmb{x}$ is positive. Hence, $\OWA_{\pmb{w}}(\pmb{v}_{\pmb{b}}(\pmb{x}))\leq 0$.
 
 Assume that $\OWA_{\pmb{w}}(\pmb{v}_{\pmb{b}}(\pmb{x}))\leq 0$ for some feasible solution $\pmb{x}\in \mathcal{X}$.  Because the costs under scenarios are nonnegative, there must be at least $r$ scenarios, say $\pmb{c}_1,\dots, \pmb{c}_r$ under which the cost of $\pmb{x}$ is~0. We get $\pmb{c}_i^{T}\pmb{x}=0$ if and only if $x_j=0$, $\overline{x}_j=1$  ($\overline{x}_j=0$, $x_j=1$) if $q_j\in \mathcal{C}_i$ ($\overline{q}_j\in \mathcal{C}_i$). It is possible for some $j\in [s]$ that neither $x_j$ nor $\overline{x}_j$ appears in the clauses $\mathcal{C}_1,\dots, \mathcal{C}_r$ corresponding to $\pmb{c}_1,\dots, \pmb{c}_r$ (so $x_j$ can be fractional). In this case, we assign any value to $q_j$ which does not change the values of $\mathcal{C}_1,\dots,\mathcal{C}_r$ (they are still not satisfied). This defines a truth assignment to $q_1,\dots, q_s$ under which at least $r$ clauses are not satisfied.
 
 The hardness of approximation follows from the fact that any $f(n)$-approximation algorithm for $\textsc{Owa}~\mathcal{P}$ could be used to verify in polynomial time if $\OWA_{\pmb{w}}(\pmb{v}_{\pmb{b}}(\pmb{x}))$ is positive.
\end{proof}

\section{Robust combinatorial problems with the OWA criterion}
\label{sec:combinatorial}

In this section, we consider the case where $\mathcal{X}\subseteq \{0,1\}^n$, so we discuss the class of combinatorial optimization problems. Observe that $\mathcal{X}$ is not a convex set and program~(\ref{mipowa1}) is only a mixed integer programming one, which is, in general, not polynomially solvable. However, the formulation~(\ref{mipowa1}) has much fewer binary variables than~(\ref{mipowa}), so the problem with non-increasing weights is still more tractable. 
It turns out that $\textsc{Owa}~\mathcal{P}$ is NP-hard for most basic combinatorial problems $\mathcal{P}$ even if $K=2$. For example, it is NP-hard for $K=2$, $\pmb{b}=\pmb{0}$, and any $\pmb{w}=(w_1,w_2)$ such that $w_1>w_2$ if $\mathcal{P}$ is the shortest path problem~\cite{KZ15}. In some cases,  $\textsc{OWA}~\mathcal{P}$ can be solved in pseudopolynomial time when $K$ is constant. Suppose that matrix~$\pmb{C}$ is integral and we can enumerate all possible vectors $\pmb{u}\in \{0,\dots,UB\}^K$, where $\pmb{u}=\pmb{v}_{\pmb{b}}(\pmb{x})$, for some solution $\pmb{x}\in \mathcal{X}$ and $UB$ is an upper bound on the components of $\pmb{v}_{\pmb{b}}(\pmb{x})$. We then find an optimal solution to $\textsc{Owa}~\mathcal{P}$ by choosing the vector $\pmb{u}$ with the minimum value of $\OWA_{\pmb{w}}(\pmb{u})$. In some cases, for example, when $\mathcal{P}$ is the shortest path problem, all vectors $\pmb{u}$ can be enumerated in pseudopolynomial time when $K$ is constant~\cite{ABV10}. Furthermore, using the reasoning from~\cite{KZ15}, the pseudopolynomial algorithm can be converted into an FPTAS under the additional assumption that the weights are non-increasing. However, the obtained algorithms are exponential in $K$, so their practical applicability is limited. In the following, we consider the case when $K$ is a part of the input.

Let us remark more on the case when the solution regrets are aggregated. If the underlying problem $\mathcal{P}$ is NP-hard, then computing the vector $\pmb{b}$ containing the optimal solution costs under scenarios is also NP-hard. Furthermore, it is easy to see that $\textsc{Owa}~\mathcal{P}$ is not approximable even if $K=1$. Indeed, solving $\mathcal{P}$ is equivalent to computing a solution $\pmb{x}$ whose regret is equal to~0. Hence, any approximation algorithm for $\textsc{Owa}~\mathcal{P}$ could be used to solve $\mathcal{P}$. We can overcome this obstacle in two ways. We can assume that $\mathcal{P}$ is polynomially solvable, or vector $\pmb{b}$ is given explicitly as a part of the input. If $\mathcal{P}$ is NP-hard, then $\pmb{b}$ can be a vector of some lower bounds on solution costs, which can be computed efficiently.

\subsection{Some hardness results}
\label{sec:hardness}

To obtain some hardness results on $\textsc{Owa}~\mathcal{P}$, we use the following known result:

\begin{thm}[\cite{KZ09, KZ11}]
\label{thmappr}
 The $\textsc{Min-max}~\mathcal{P}$ problem 
  is strongly NP-hard and
 hard to approximate within $O(\log^{1-\epsilon}K)$ for any $\epsilon>0$,
 unless NP\,$\subseteq \text{DTIME}(n^{\text{polylog}(n)})$,
 when $\mathcal{P}$ is the shortest path, minimum spanning tree,  minimum assignment, or 
 minimum s-t cut. 
\end{thm}
Theorem~\ref{thmappr} remains true when $\pmb{b}=\pmb{0}$.
To extend the hardness results for  $\textsc{Owa}~\mathcal{P}$, we first prove the following proposition:
\begin{prop}
\label{propha}
If  $\textsc{Owa}~\mathcal{P}$ with non-increasing weights is approximable within $\gamma\geq 1$, then $\textsc{Min-max}~\mathcal{P}$ is approximable within $\frac{\gamma}{w_1}$.
\end{prop}
\begin{proof}
	Let $\pmb{x}^*$ be an optimal solution to $\textsc{Min-max}~\mathcal{P}$, $\pmb{x}'$ and $\pmb{x}''$ be an optimal and $\gamma$-approximate solution to~$\textsc{Owa}~\mathcal{P}$, respectively. By 
	Proposition~\ref{eglp} (set $p=\infty$ and $q=1$, since $\sum_{i\in[K]}w_i=1$,  $\lVert \pmb{w}\rVert_{1}=1$ and,
	in consequence, $\rho=\frac{1}{w_1}$) we have
	$$ \lVert \pmb{v}_{\pmb{b}}(\pmb{x}^*) \rVert _{\infty}\geq \OWA_{\pmb{w}}(\pmb{v}_{\pmb{b}}(\pmb{x}^*))\geq \OWA_{\pmb{w}}(\pmb{v}_{\pmb{b}}(\pmb{x}')) \geq \frac{1}{\gamma}\OWA_{\pmb{w}}(\pmb{v}_{\pmb{b}}(\pmb{x}''))\geq
	\frac{w_1}{\gamma} \lVert \pmb{v}_{\pmb{b}}(\pmb{x}'') \rVert _{\infty}.
	 $$
Hence $\pmb{x}''$ is a $\frac{\gamma}{w_1}$-approximate solution to $\textsc{Min-max}~\mathcal{P}$.
\end{proof}
Let $\pmb{w}$ be a vector of non-increasing weights. We can characterize $\pmb{w}$ by considering the largest weight $w_1$,
obviously,
$w_1 \in [\frac{1}{K},1]$.
This weight can be either a constant or a function of $K$.  The following corollaries  are a direct   consequence of Proposition~\ref{propha} and Theorem~\ref{thmappr}:
\begin{cor}
	If $w_1$ is a constant, then the  $\textsc{Owa}~\mathcal{P}$ problem with non-increasing weights 
	 is strongly NP-hard
	 and
	 hard to approximate within $O(\log^{1-\epsilon}K)$ for any $\epsilon>0$, unless NP\,$\subseteq \text{DTIME}(n^{\text{polylog}(n)})$,
	 when $\mathcal{P}$ is the shortest path, minimum spanning tree,  minimum assignment,
	 or 
 minimum s-t cut. 
 \label{cw1const}
\end{cor}
\begin{cor}
\label{cor2}
 If $w_1\geq \frac{1}{\log^{1-\epsilon} K}$ for some $\epsilon>0$, then the $\textsc{Owa}~\mathcal{P}$ problem with non-increasing weights
   is strongly NP-hard and
  hard to approximate within any constant factor, unless NP\,$\subseteq \text{DTIME}(n^{\text{polylog}(n)})$, when~$\mathcal{P}$ is the shortest path, minimum spanning tree, minimum assignment,
 or 
 minimum s-t cut. 
 \end{cor}
 \begin{proof}
 	Assume that there is a $\gamma$-approximation algorithm for $\textsc{Owa}~\mathcal{P}$ for some constant $\gamma\geq 1$. Then, according to Proposition~\ref{propha}, there is a $O(\log^{1-\epsilon} K)$-approximation algorithm for $\textsc{Min-max}~\mathcal{P}$, for some $\epsilon>0$, which contradicts Theorem~\ref{thmappr}.
 \end{proof}
 It has been shown in~\cite{KZ15} that $\textsc{Owa}~\mathcal{P}$ with non-increasing weights and $\pmb{b}=\pmb{0}$ is approximable within $w_1K$ when $\mathcal{P}$ is polynomially solvable.  Therefore, $\textsc{Owa}~\mathcal{P}$ is then approximable within a constant factor  if $w_1=\Theta(\frac{1}{K})$. On the other hand, by Corollary~\ref{cor2}, the problem is hard to approximate within any constant factor if $w_1=\Omega(\frac{1}{\log^{1-\epsilon} K})$ for some  $\epsilon>0$ or
 by Corollary~\ref{cw1const}, the problem is hard to approximate within
$O(\log^{1-\epsilon}K)$ for any $\epsilon>0$  if $w_1$ is  a constant.

\subsection{Approximation algorithms based on $p$-norm minimization}
\label{sec:approx}

In this section, we present new approximation results for $\textsc{Owa}~\mathcal{P}$ when the vector of weights $\pmb{w}$ is non-increasing. Consider first  the following auxiliary problem:
\begin{equation}
      \text{$p$-\textsc{norm}}~\mathcal{P}:\; \min_{\pmb{x}\in\mathcal{X}}\lVert\pmb{v}_{\pmb{b}}(\pmb{x})\rVert_p=\min_{\pmb{x}\in\mathcal{X}}\lVert\pmb{C}\pmb{x}-\pmb{b}\rVert_p,
      \label{lpp}
\end{equation}
where $p\in \overline{\mathbb{R}}_{+}$. For $\pmb{b}=\pmb{0}$, we get the $p$-norm minimization problem discussed, for instance, in~\cite{BCFM17}. For arbitrary $\pmb{b}$, we minimize a distance to 
 a \emph{reference point} or \emph{ideal point} (if $\pmb{b}$ is a vector of the optimal costs under scenarios). Such a problem is commonly used in multi-objective optimization, and some results in this area for combinatorial problems can be found in~\cite{BGMS17}.
 Our new approximation results are based on the following proposition:
\begin{prop}
\label{propappr}
	If $\text{$p$-\textsc{norm}}~\mathcal{P}$ is approximable within $\gamma\geq 1$, then $\textsc{Owa}~\mathcal{P}$ with non-increasing weights is approximable within $\gamma \rho$, where $\rho=K^{\frac{1}{p}}\frac{ \lVert \pmb{w}\rVert_{q} }{\lVert \pmb{w}\rVert_p}$ and $\frac{1}{p}+\frac{1}{q}=1$.
\end{prop}
\begin{proof}
	Let $\pmb{x}^*$ be an optimal solution to $\textsc{Owa}~\mathcal{P}$ and let $\pmb{x}'$ be a $\gamma$-approximate solution to $\text{$p$-\textsc{norm}}~\mathcal{P}$. Using~Proposition~\ref{eglp} we get
$$\OWA_{\pmb{w}}(\pmb{v}_{\pmb{b}}(\pmb{x}^*))\geq \frac{1}{\rho}\lVert \pmb{w}\rVert_q  \lVert \pmb{v}_{\pmb{b}}(\pmb{x}^*)\rVert_p \geq \frac{1}{\rho\gamma} \lVert \pmb{w}\rVert_q  \lVert \pmb{v}_{\pmb{b}}(\pmb{x}')\rVert_p\geq \frac{1}{\rho\gamma} \OWA_{\pmb{w}}(\pmb{v}_{\pmb{b}}(\pmb{x}'))$$
and $\pmb{x}'$ is a $\gamma\rho$-approximate solution to $\textsc{Owa}~\mathcal{P}$.
\end{proof}

Let us analyze the factor $\rho$ from Proposition~\ref{propappr}. Table~\ref{tab:rho} shows the values of $\rho$ for $p\in \{1,2,\infty\}$.
If $p\le 2$, then $q=\frac{p}{p-1} \ge p$. Using Proposition~\ref{norinq}, we find
$\frac{\lVert \pmb{w} \rVert_q}{\lVert \pmb{w} \rVert_p} \le 1$ and  $\rho \le K^{1/p}$.
If $p \ge 2$, then $q \le p$. Again, using Proposition~\ref{norinq}, we obtain  $\frac{\lVert \pmb{w} \rVert_q}{\lVert \pmb{w} \rVert_p} \le K^{\frac{1}{q}-\frac{1}{p}}$ and $\rho \le K^{1/q}$. Therefore, the best upper bound from this case distinction is reached when $p=2$.  Notice that for $p=2$, the part of $\rho$ that depends of $\pmb{w}$ cancels, and $\rho=\sqrt{K}$. The value of $\rho$ can depend on the weight distribution in $\pmb{w}$. Let us recall that for non-increasing weights, we have $w_1\in [\frac{1}{K},1]$. When $w_1$ is close to $\frac{1}{K}$, then we should fix $p=1$. On the other hand, when $w_1$ is close to~1, we should choose $p=\infty$. 
\begin{table}[htb]
\caption{Value of $\rho$ from Proposition~\ref{propappr} for different values of $p$.}\label{tab:rho}
\begin{center}
\begin{tabular}{r|ccc}
$p$ & 1 & 2 & $\infty$ \\
\hline
$\rho$ & $w_1K$ & $\sqrt{K}$  & $\frac{1}{w_1}$
\end{tabular}
\end{center}
\end{table}

\begin{figure}[ht]
	\centering
	\includegraphics[height=6cm]{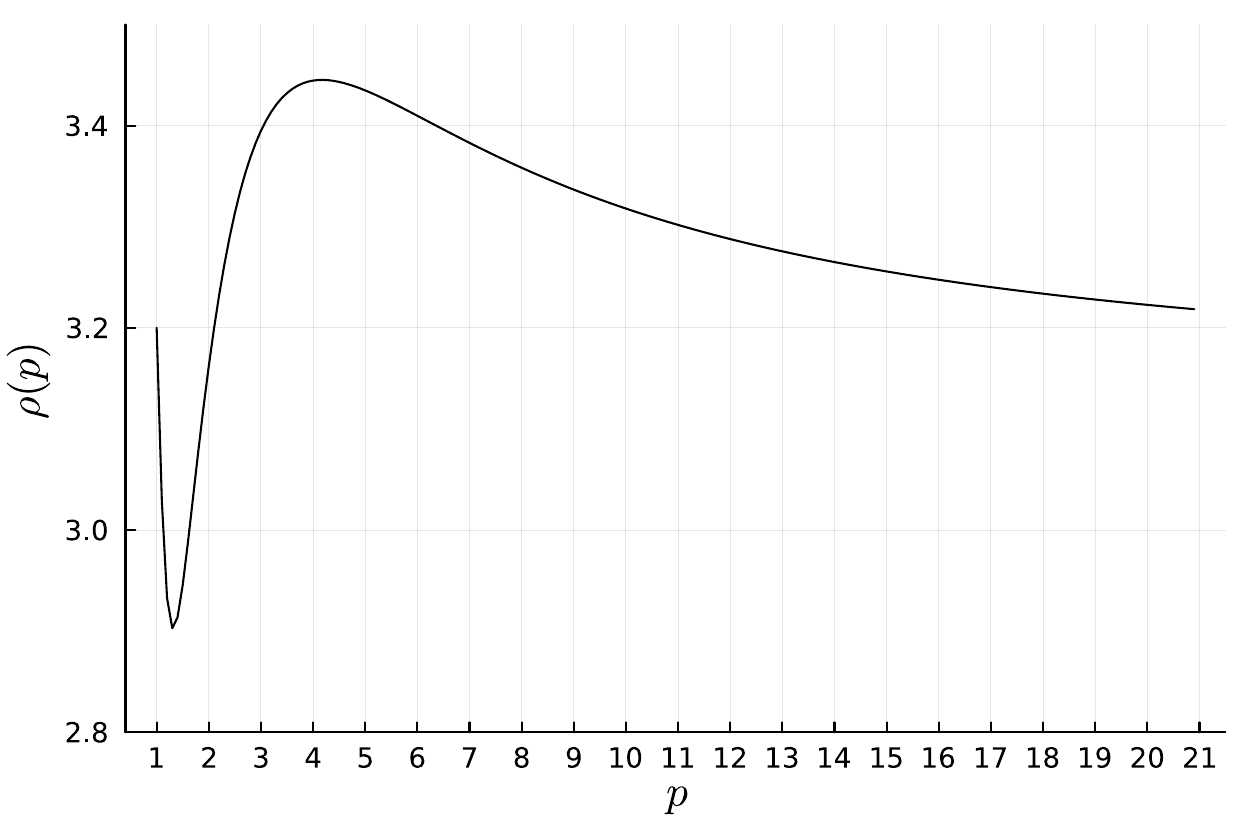} 
	\caption{The function $\rho(p)$ for $\pmb{w}=(0.32, 0.22, 0.12, 0.12, 0.12, 0.1, 0,0,0,0)$.} \label{fig2}
\end{figure}

  Using a more careful analysis, we should compute $p^*$ by solving the problem.

\begin{equation}
\label{bestrho}
p^*\in \mathrm{arg}\,\inf_{p\in [1,\infty)} \rho(p)= K^{\frac{1}{p}}\frac{ \lVert \pmb{w}\rVert_{\frac{p}{p-1}} }{\lVert \pmb{w}\rVert_p}.
\end{equation}
It is worth noting that the value of $p^*$ depends on the whole weight distribution in $\pmb{w}$. To see that $\rho(p^*)$ can be smaller than $\min\{\rho(1), \rho(2),  \rho(\infty)\}$ let us consider
\[
\pmb{w}=(0.32, 0.22, 0.12, 0.12, 0.12, 0.1, 0,0,0,0)
\] with $K=10$. The function $\rho(p)$ for $p\geq 1$ is shown in Figure~\ref{fig2}. The function $\rho(p)$ tends to $\frac{1}{w_1}=3.125$ as $p\rightarrow \infty$.
Thus $p^*\approx 1.3$ and $\rho(p^*)\approx2.9$, while $\rho(1)=3.2$, $\rho(2)=3.16$ and $\rho(\infty)=3.125$.

A more detailed analysis of the approximation ratio from Proposition~\ref{propappr} should take into account the value of $\gamma$. In general, $\gamma$ may also depend on $p$. Therefore, for a particular problem we should find $p^*$ minimizing the product $\gamma(p)\rho(p)$, where $\gamma(p)$ is an approximation factor of the problem $\text{$p$-\textsc{norm}}~\mathcal{P}$. Let us analyze some special cases of $p$ in more detail.

\begin{prop}
\label{prop1norm}
 If one can approximate $\text{1-\textsc{norm}}~\mathcal{P}$ within a factor of $\alpha\geq 1$, then $\textsc{Owa}~\mathcal{P}$ with non-increasing weights is approximable within $\alpha w_1 K$.
 \end{prop}
 Observe that   $\text{1-\textsc{norm}}~\mathcal{P}$ is polynomially solvable if $\mathcal{P}$ can be solved in polynomial time. Indeed, by the assumption that $\pmb{C}\pmb{x}-\pmb{b}\geq \pmb{0}$ for each $\pmb{x}$, we get $\lVert\pmb{v}_{\pmb{b}}(\pmb{x})\rVert_1=\sum_{i\in [K]} (\pmb{c}_i^T\pmb{x}-b_i)=\hat{\pmb{c}}^T\pmb{x}+\hat{b}$, where $\hat{\pmb{c}}=\sum_{i\in [K]} \pmb{c}_i$ and $\hat{b}=\sum_{i\in [K]} b_i$. Thus, it is enough to solve the problem $\mathcal{P}$ for the cost vector~$\hat{\pmb{c}}$. Therefore, if $\mathcal{P}$ is polynomially solvable, then $\alpha=1$.
 A result similar to Proposition~\ref{prop1norm} has been shown in~\cite{KZ15}.
 Since $w_1\in [\frac{1}{K},1]$, the approximation ratio can be $O(K)$ when $\alpha$ is constant. A better approximation ratio can be achieved when we have a $\beta$-approximation algorithm for the $\text{2-\textsc{norm}}~\mathcal{P}$ problem. Fixing $p=2$ in Proposition~\ref{propappr}, we conclude that $\rho=\sqrt{K}$ and $\textsc{Owa}~\mathcal{P}$ is then approximable within $\beta \sqrt{K}$ which is $O(\sqrt{K})$ if $\beta$ is constant. 
 \begin{prop}
 	\label{prop2norm}
 If one can approximate $\text{2-\textsc{norm}}~\mathcal{P}$ within a factor of $\beta\geq 1$, then $\textsc{Owa}~\mathcal{P}$ with non-increasing weights is approximable within $\beta \sqrt{K}$.
 \end{prop}
 Finally, having a $\gamma$-approximation algorithm for  $\text{$\infty$-\textsc{norm}}~\mathcal{P}$ gives us $\rho=\frac{1}{w_1}$, which yields a $\frac{\gamma}{w_1}$-approximation algorithm for $\textsc{Owa}~\mathcal{P}$. 
  \begin{prop}
 	\label{propinftynorm}
 If one can approximate $\text{$\infty$-\textsc{norm}}~\mathcal{P}$ within a factor of $\gamma\geq 1$, then $\textsc{Owa}~\mathcal{P}$ with non-increasing weights is approximable within $\frac{\gamma}{w_1}$.
 \end{prop}
 Combining the three approximation algorithms, we obtain the following characterization of the approximability of $\textsc{Owa}~\mathcal{P}$ that depends on $w_1\in [\frac{1}{K},1]$ and $K$:
\begin{cor}
\label{propappr1}
	If  $\text{p-\textsc{norm}}~\mathcal{P}$ for $p=1,2,\infty$ is approximable within $\alpha$, $\beta$, $\gamma$, respectively, then $\textsc{Owa}~\mathcal{P}$ with non-increasing weights is approximable within $\min\{\alpha w_1 K, \beta \sqrt{K}, \frac{\gamma}{w_1}\}$.
\end{cor}
For many basic combinatorial problems, $\alpha$ and $\beta$ are constant (in particular, $\alpha=1$ if $\mathcal{P}$ is polynomially solvable). However, $\gamma$ typically is not constant, as $\text{$\infty$-\textsc{norm}}~\mathcal{P}$ is equivalent to $\textsc{Min-max}~\mathcal{P}$  (see Theorem~\ref{thmappr}). In general, if $\mathcal{P}$ can be solved in polynomial time, then $\textsc{Min-max}~\mathcal{P}$ can be approximated within $\gamma=K$ (see~\cite{ABV06}). A better approximation ratio can be achieved for particular problems, and we analyze such cases in the next sections.

\subsection{Application to the matroidal, the shortest path, and the minimum Steiner tree problems}
\label{sec:appl}

In this section, we apply the results from Sections~\ref{sec:hardness} and~\ref{sec:approx} to some particular problems. We first consider the case with $\pmb{b}=\pmb{0}$, so only the solution costs over scenarios are aggregated. Next, we consider the more general case with $\pmb{b}>\pmb{0}$. We use the following known results:
\begin{thm}[\cite{BCFM17}]
\label{tacop}
There exist algorithms that approximate the $\text{$p$-\textsc{norm}}~\mathcal{P}$ problem for $\pmb{b}=\pmb{0}$ within a factor of
$O(\min\{p,\log K\})$, for $p\in \overline{\Rset}_{+}$, if $\mathcal{X}$ describes 
the sets of feasible solutions of 
matroidal problems, the shortest path problem, or the minimum Steiner tree problem.
\end{thm} 
Observe that $\textsc{Min-max}~\mathcal{P}$ (i.e. the case with $p=\infty$), is then approximable within $O(\log K)$, so $\gamma$ from Corollary~\ref{propappr1} is not a constant (it depends on $K$). 

\subsubsection{The case of $\pmb{b}=\pmb{0}$}

Theorem~\ref{tacop}, together with Proposition~\ref{propappr} imply the following result:
\begin{thm} \label{thmowa}
	For each $p\in \overline{\Rset}_{+}$, the $\textsc{Owa}~\mathcal{P}$ problem with non-increasing weights and $\pmb{b}=\pmb{0}$ is approximable within $O\left(K^{\frac{1}{p}}\frac{ \lVert \pmb{w}\rVert_{\frac{p}{p-1}} }{\lVert \pmb{w}\rVert_p}\min\{p,\log K\}\right)$ when $\mathcal{P}$ is a matroidal problem, the shortest path problem, or the minimum Steiner tree problem.
\end{thm}
By considering the special cases $p=1,2,\infty$ (see also Corollary~\ref{propappr1}), we get
\begin{cor}
\label{corapprspt}
The $\textsc{Owa}~\mathcal{P}$ problem with non-increasing weights and $\pmb{b}=\pmb{0}$ is approximable within  $O(\min\{w_1K, \sqrt{K}, \frac{\log K}{w_1}\})$ when $\mathcal{P}$ is a matroidal problem, the shortest path problem, or the minimum Steiner tree problem.
\end{cor}

The approximation ratio from Corollary~\ref{corapprspt}   is $O(\sqrt{K})$, which significantly improves the known $w_1 K$-approximation algorithm known to date~\cite{KZ15}.
Corollary~\ref{corapprspt}, together with Theorem~\ref{thmappr}, allows us to provide the following characterization of the approximability of $\textsc{OWA}~\mathcal{P}$ for the considered problems $\mathcal{P}$. If $w_1$ is a constant, then the problem is approximable within $O(\log K)$ but hard to approximate within $O(\log^{1-\epsilon} K)$ for any $\epsilon>0$. If $w_1=\Theta(\frac{1}{K^{1/q}})$, for some $q\geq 1$, then the problem is approximable within $O(\min\{K^{\frac{q-1}{q}},\sqrt{K}, K^{\frac{1}{q}}\cdot\log K \})$. In particular, if $q=1$, the problem is approximable within a constant factor~$O(1)$. Asymptotically, the approximation ratio attains the largest value for $q=2$, which is $O(\sqrt{K})$. Therefore, the worst weight distribution for the approximation algorithm is when 
$w_1=\Theta(\frac{1}{\sqrt{K}})$. Finally, if $w_1=\Omega(\frac{1}{\log^{1-\epsilon} K})$ for some constant~$\epsilon>0$, then the problem is approximable within $O(\log^{2-\epsilon} K)$. However, in this case, it is hard to approximate within any constant factor (see Corollary~\ref{cor2}). Notice that there is still room for improvement for the approximation results when $w_1$ depends on $K$.

\begin{figure}[ht]
	\centering
	\includegraphics[height=6cm]{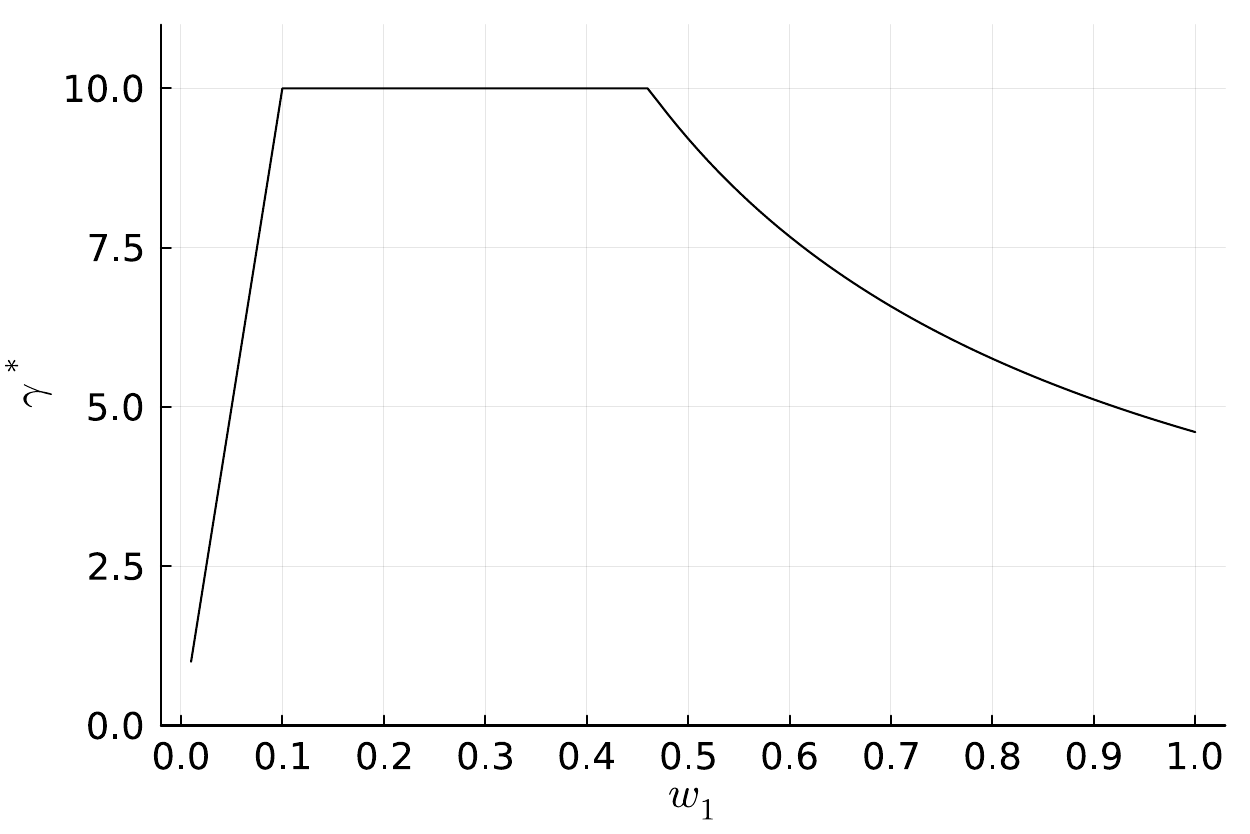}
	\caption{The ratio $\min\{w_1K, \sqrt{K}, \frac{\log K}{w_1}\}$ for $K=100$, depending on $w_1$.} \label{fig1}
\end{figure}

Figure~\ref{fig1} shows the value of $\min\{w_1K, \sqrt{K}, \frac{\log K}{w_1}\}$ for $K=100$ and $w_1\in [\frac{1}{100},1]$. Notice that for small $w_1$, i.e. when the distribution of the weights is close to uniform, we should use $p=1$. If $w_1$ is large, then we should use $p=\infty$. Finally, for intermediate $w_1$, the best ratio is achieved using $p=2$.

\subsubsection{The case of  $\new{\pmb{b}>\pmb{0}}$}

Unfortunately, the approximation results obtained in~\cite{BCFM17} cannot directly be applied when~$\new{\pmb{b}>\pmb{0}}$. To apply 
them for this more general case, we need the following proposition:
\begin{prop}
\label{propapprb}
	Assume that $\lVert \pmb{C}\pmb{x}-\pmb{b}\rVert_{p}\geq 1$ for each $\pmb{x}\in \mathcal{X}$.
	Let $\pmb{x}^*\in \mathcal{X}$ minimize $
	\Vert\pmb{C}\pmb{x}\rVert_{p}$ for $p\in \overline{\Rset}_{+}$. Then for each $\pmb{x}\in \mathcal{X}$.
	$$\lVert\pmb{C}\pmb{x}^*-\pmb{b}\rVert_{p}\leq (\lVert\pmb{b}\rVert_{p}+1)\lVert \pmb{C}\pmb{x}-\pmb{b}\rVert_{p}.$$
\end{prop}
\begin{proof}
	For each solution $\pmb{x}\in \mathcal{X}$, we get $\lVert\pmb{C}\pmb{x}^*-\pmb{b}\rVert_{p}\leq\lVert\pmb{C}\pmb{x}^*\rVert_{p}\leq \lVert\pmb{C}\pmb{x}\rVert_{p}\leq \lVert\pmb{C}\pmb{x}-\pmb{b}\rVert_{p}+\lVert\pmb{b}\rVert_{p}\leq (\lVert\pmb{b}\rVert_{p}+1)\lVert \pmb{C}\pmb{x}-\pmb{b}\rVert_{p}$, where the first  inequality follows from the fact that 
	 $\pmb{c}_i^T\pmb{x}\geq b_i$, $i\in [K]$, for every $\pmb{x}\in \mathcal{X}$,
	the third inequality results from the triangle inequality and the assumption that $\lVert \pmb{C}\pmb{x}-\pmb{b}\rVert_{p}\geq 1$
	gives 
	the last inequality.
\end{proof}

\begin{thm}
\label{thmapprb}
Assume that $\lVert \pmb{C}\pmb{x}-\pmb{b}\rVert_{p}\geq 1$ for each $\pmb{x}\in \mathcal{X}$.
For each $p\in \overline{\Rset}_{+}$, the $\textsc{Owa}~\mathcal{P}$ problem with non-increasing weights is approximable within $O\left(K^{\frac{1}{p}}\frac{ \lVert \pmb{w}\rVert_{\frac{p}{p-1}} }{\lVert \pmb{w}\rVert_p}(\lVert \pmb{b} \rVert_{p}+1)\min\{p,\log K\}\right)$ when $\mathcal{P}$ is a matroidal problem, the shortest path problem, or the minimum Steiner tree problem.
\end{thm}
\begin{proof}
	  It follows from Theorem~\ref{tacop} that we can obtain an $O(\min\{p, \log K\})$-approximate solution to the problem of minimizing $\lVert \pmb{C}\pmb{x}\rVert_{p}$. By Proposition~\ref{propapprb},  this solution is $O((\lVert \pmb{b} \rVert_{p}+1)\min\{p, \log K\})$-approximate to the problem of minimizing $\lVert \pmb{C}\pmb{x}-\pmb{b}\rVert_{p}$. Hence and by Proposition~\ref{propappr}
	we have the assertion of the theorem.
\end{proof}
Assume that all the components of $\pmb{C}$ and $\pmb{b}$ are integral. Then $\lVert \pmb{C}\pmb{x}-\pmb{b}\rVert_{p}< 1$  if and only if $\pmb{c}_i^T\pmb{x}-b_i=0$ for each $i\in [K]$.  This is equivalent to $\sum_{i\in [K]} (\pmb{c}_i^T\pmb{x}-b_i)=0$, because $\pmb{c}_i^T\pmb{x}-b_i\geq 0$ for each $i\in [K]$. Therefore, if the underlying problem $\mathcal{P}$ is polynomially solvable, then the assumption of Theorem~\ref{thmapprb} can be checked in polynomial time by finding an optimal solution $\pmb{x}^*$ to $\mathcal{P}$ for the cost vector $\hat{\pmb{c}}=\sum_{i\in [K]}\pmb{c}_i$ and checking if $\sum_{i\in [K]} (\pmb{c}_i^T\pmb{x}^*-b_i)=0$. Notice also that $\pmb{x}^*$ is then an optimal solution to $\textsc{Owa}~\mathcal{P}$ because the value of the objective in this problem is nonnegative. Hence, for integral data, the assumption of Theorem~\ref{thmapprb} is not very restrictive.

We now show an additional application of Theorem~\ref{thmapprb}.
Let us construct instance $(\hat{\pmb{C}},\hat{\pmb{b}})$ by dividing each component of $\pmb{C}$ and $\pmb{b}$ by 
$\lVert \pmb{b}\rVert_{p}>0$. 
It is easy to see that $\pmb{c}_i^T\pmb{x}-b_i=\lVert \pmb{b}\rVert_{p}(\hat{\pmb{c}}_i^T\pmb{x}-\hat{b}_i)$ for each $i\in [K]$ and thus $\hat{\pmb{c}}_i^T\pmb{x}-\hat{b}_i\geq 0$, $i\in [K]$, for each feasible solution $\pmb{x}\in \mathcal{X}$.
Clearly $\lVert \hat{\pmb{b}}\rVert_{p}=1$.
\begin{prop}
	If $\textsc{Owa}~\mathcal{P}$ for the instance $(\hat{\pmb{C}},\hat{\pmb{b}})$ is approximable within $\gamma$, then also the problem for the instance $(\pmb{C},\pmb{b})$ is approximable within $\gamma$.
	\label{pscal}
\end{prop}
\begin{proof}
	Let $\OWA_{\pmb{w}}(\hat{\pmb{v}}_{\hat{\pmb{b}}}(\pmb{x}))$ be the value of OWA of $\pmb{x}$ for the scaled instance  $(\hat{\pmb{C}},\hat{\pmb{b}})$.
	Thus $\OWA_{\pmb{w}}(\pmb{v}_{\pmb{b}}(\pmb{x}))=\lVert \hat{\pmb{b}}\rVert_{p} \OWA_{\pmb{w}}(\hat{\pmb{v}}_{\hat{\pmb{b}}}(\pmb{x}))$. Let $\pmb{x}^*$ be a $\gamma$-approximate solution for $(\hat{\pmb{C}},\hat{\pmb{b}})$. Then
	for each $\pmb{x}\in \mathcal{X}$
	$$\OWA_{\pmb{w}}(\pmb{v}_{\pmb{b}}(\pmb{x}^*))=\lVert \hat{\pmb{b}}\rVert_{p}\OWA_{\pmb{w}}(\hat{\pmb{v}}_{\hat{\pmb{b}}}(\pmb{x}^*))\leq \gamma \lVert \hat{\pmb{b}}\rVert_{p} \OWA_{\pmb{w}}(\hat{\pmb{v}}_{\hat{\pmb{b}}}(\pmb{x}))\leq \gamma \OWA_{\pmb{w}}(\pmb{v}_{\pmb{b}}(\pmb{x}))$$
and the proposition follows.
\end{proof}

We now apply Theorem~\ref{thmapprb} to the scaled instance $(\hat{\pmb{C}},\hat{\pmb{b}})$.  We have to check first the assumption that  $\lVert \hat{\pmb{C}}\pmb{x}-\hat{\pmb{b}}\rVert_{p}\geq 1$ for each $\pmb{x}\in \mathcal{X}$. This is not an easy task in general (notice that now the components of  $(\hat{\pmb{C}},\hat{\pmb{b}})$ need not be integral).
Since  $\hat{\pmb{c}}_i^T\pmb{x}\geq \hat{b}_i$, $i\in [K]$, for every $\pmb{x}\in \mathcal{X}$ and 
$\lVert \hat{\pmb{b}}\rVert_{p}=1$,
 \[
 \lVert \hat{\pmb{C}}\pmb{x}-\hat{\pmb{b}}\rVert_{p}\geq \lVert \hat{\pmb{C}}\pmb{x}\rVert_{p}- \lVert \hat{\pmb{b}}\rVert_{p}=
 \lVert \hat{\pmb{C}}\pmb{x}\rVert_{p}-1\geq  \lVert \hat{\pmb{C}}\hat{\pmb{x}}\rVert_{p}-1\geq
 \frac{1}{\eta} \lVert \hat{\pmb{C}}\pmb{x}'\rVert_{p}-1
 \geq 0,
 \]
 where  $\hat{\pmb{x}}\in \mathcal{X}$ is an optimal solution to the problem of minimizing $\lVert \hat{\pmb{C}}\pmb{x}\rVert_{p}$
 and  $\pmb{x}'\in \mathcal{X}$ is its $\eta$-approximate solution to this problem.
 Accordingly, we  only need to find~$\pmb{x}'$ by using a $\eta$-approximation algorithm (by Theorem~\ref{tacop}, 
  $\eta=O(\min\{p,\log K\})$) and check if $ \frac{1}{\eta} \lVert \hat{\pmb{C}}\pmb{x}'\rVert_{p}\geq 1$ holds.
  If this inequality is satisfied, then by Proposition~\ref{pscal} and Theorem~\ref{thmapprb}, $\pmb{x}'$ is an
  $O\left(K^{\frac{1}{p}}\frac{ \lVert \pmb{w}\rVert_{\frac{p}{p-1}} }{\lVert \pmb{w}\rVert_p}\eta\right)$-
  approximate solution for the original problem.
  Observe that
  setting $p=2$ in Theorem~\ref{thmapprb} leads to a $O(\sqrt{K})$-approximation algorithm for the considered problems.
  Another method for checking if the assumption of Theorem~\ref{thmapprb} is met
 consists in solving some computationally efficient relaxation of the $\text{$p$-\textsc{norm}}~\mathcal{P}$ problem. For instance, a convex relaxation by simply replacing binary constraints $x_i\in \{0,1\}$ with $0\leq x_i\leq 1$ for $i\in [n]$ in the description of $\mathcal{X}$.  We get a lower bound $LB$ on  the optimal objective value of $\text{$p$-\textsc{norm}}~\mathcal{P}$. 
  Obviously,
  the assumption is satisfied if $LB\geq 1$.

\subsection{Scenario aggregation}
\label{sec:aggregate}

In this section, we recall another approach to approximate $\textsc{Owa}~\mathcal{P}$ for non-increasing weights. The idea (see~\cite{CGKZ20}) consists in reducing the number of scenarios and solving a smaller instance using the formulations~(\ref{mipowa1}) or~(\ref{mipowa}). Let $K$ be a multiple of $\ell$. Define $\overline{\pmb{C}} \in \Rset^{K/\ell\times n}_{+}$,  $\overline{\pmb{b}}\in \Rset_{+}^{K/\ell}$,  and  $\overline{\pmb{w}}\in  \Rset_{+}^{K/\ell}$, where $\overline{\pmb{c}}_i= \sum_{j\in[\ell]} \pmb{c}_{(i-1)\ell  + j}/\ell$, $\overline{b}_i=\sum_{j\in[\ell]} b_{(i-1)\ell  + j}/\ell$, and $\overline{w}_i = w_{(i-1)\ell+1} + \ldots + w_{i\ell}$.  The instance $(\overline{\pmb{C}},\overline{\pmb{b}},\overline{\pmb{w}})$ is an aggregated instance 
$(\pmb{C},\pmb{b},\pmb{w})$ of $\textsc{Owa}~\mathcal{P}$.
Let
$$\overline{\pmb{v}}_{\overline{\pmb{b}}}(\pmb{x})=(\overline{\pmb{c}}_1^T\pmb{x}-\overline{b}_1,\dots, \overline{\pmb{c}}_{K/\ell}^T\pmb{x}-\overline{b}_{K/\ell}).$$
It is easy to see that $\overline{\pmb{c}}_i^T\pmb{x}-\overline{b}_i =\sum_{j\in[\ell]} (\pmb{c}_{(i-1)\ell  + j}^T\pmb{x}-b_{(i-1)\ell+j})/\ell \geq 0$ for each $\pmb{x}\in \mathcal{X}$. Therefore, using Proposition~\ref{lemaggr}, we get the following result:
\begin{prop}
\label{propaggr1}
If the vector of weights $\pmb{w}$ is non-increasing, then for each solution $\pmb{x}\in \mathcal{X}$ the inequalities
\begin{equation}
\OWA_{\overline{\pmb{w}}}(\overline{\pmb{v}}_{\overline{\pmb{b}}}(\pmb{x})) \le \OWA_{\pmb{w}}(\pmb{v}_{\pmb{b}}(\pmb{x})) \le \ell\phi\OWA_{\overline{\pmb{w}}}(\overline{\pmb{v}}_{\overline{\pmb{b}}}(\pmb{x}))
\end{equation}
hold, where $\phi = \max_{k\in[K/\ell]} \left(\sum_{i\in[k]} w_i / \sum_{i\in[k]} \overline{w}_i\right)$.
\end{prop}

\begin{cor}
\label{corappr}
	Assume that the vector of weights $\pmb{w}$ is non-increasing.
	If $\pmb{x}'$ is an optimal solution to $\textsc{Owa}~\mathcal{P}$ for the aggregated instance $(\overline{\pmb{C}},\overline{\pmb{b}},\overline{\pmb{w}})$, then $\pmb{x}'$ is an $\ell\phi$-approximate solution to $\textsc{Owa}~\mathcal{P}$ for the instance $(\pmb{C},\pmb{b},\pmb{w})$, where $\phi = \max_{k\in[K/\ell]} \left(\sum_{i\in[k]} w_i / \sum_{i\in[k]} \overline{w}_i\right)$.
\end{cor}
\begin{proof}
Let $\pmb{x}^*\in \mathcal{X}$ be an optimal solution to $\textsc{Owa}~\mathcal{P}$ for $(\pmb{C},\pmb{b},\pmb{w})$. Then, using Proposition~\ref{propaggr1} we get
$$\OWA_{\pmb{w}}(\pmb{v}_{\pmb{b}}(\pmb{x}^*))\geq \OWA_{\overline{\pmb{w}}}(\overline{\pmb{v}}_{\overline{\pmb{b}}}(\pmb{x}^*))\geq \OWA_{\overline{\pmb{w}}}(\overline{\pmb{v}}_{\overline{\pmb{b}}}(\pmb{x}'))\geq \frac{1}{\ell \phi} \OWA_{\pmb{w}}(\pmb{v}_{\pmb{b}}(\pmb{x}'))$$
and the corollary follows.
\end{proof}
Using Corollary~\ref{corappr}, we can significantly reduce the size of the problem, preserving some approximation guarantee. For, example, when $\ell=2$, we can reduce the size of $(\pmb{C}, \pmb{b}, \pmb{w})$ by 50\%. Solving the reduced (aggregated) instance $(\overline{\pmb{C}},\overline{\pmb{b}},\overline{\pmb{w}})$ we get a $2\phi$-approximate solution to the original problem. Notice that $\phi\leq 1$, so the approximation ratio in this case is at most~2. A more detailed analysis of the aggregation method and its experimental evaluation can be found in~\cite{CGKZ20}.

\section{Experiments}
\label{sec:experiments}

In this section, our attention is directed toward the distinct case, where $\pmb{b}$ represents the vector of optimal objective values in each scenario. For the purpose of these experiments, we refer to this setting as Ordered Weighted Averaging Regret (OWAR). We describe three types of computational experiments.

In the first experiment, we examine the theoretical results derived from Theorem~\ref{thmowa} (pertaining to OWA) and Theorem~\ref{thmapprb} (related to OWAR). We employ the greedy algorithm that starts with an empty solution, and in each iteration adds an item to the solution that minimizes the current $p$-norm. The solution of this algorithm is then evaluated in the context of OWA and in OWAR. We validate our theoretical findings and conduct a comparative analysis by assessing the results under various $p$-norms and weights. 
This investigation aims to provide practical insights into the performance of the proposed frameworks in varied decision-making scenarios.

In the second experiment, we consider the advantage that OWAR offers as a generalized decision-making criterion by comparing the solutions we find to those of min-max regret and classic OWA.

In the third experiment, our focus shifts to scenario aggregation algorithms. We evaluate the theoretical results derived from Corollary \ref{corappr} by implementing the $\ell$-aggregation and $K$-means algorithms. This examination seeks to elucidate the effectiveness of these algorithms in aggregating scenarios and their implications for decision-making. The subsequent sections will detail the experimental setup and results, providing a comprehensive understanding of the diverse facets explored in our study.
All experiments were performed on a machine with a 6-core Intel i7 2.6 GHz processor. The implementation was done in Python 3.11, utilizing various libraries including Gurobi 10.0, NumPy, and scikit-learn.

\subsection{Experiment 1: Performance of greedy algorithm}
\subsubsection{Setup}\label{sec:setup1}

We consider randomly generated selection problems with $\X = \{ \pmb{x}\in\{0,1\}^n : \sum_{i\in[n]} x_i = q\}$,
which are matroidal problems.
 The approximation guarantee mentioned in Theorem~\ref{tacop} stems from a greedy algorithm, where items are packed sequentially so that the $p$-norm of the current objective vector is minimized. The guarantee $O(\min\{p,\log K\})$ is achieved by running this greedy algorithm twice; for the sake of simplicity, we run the greedy method only once, which gives a $O(p)$-guarantee. The goal is to evaluate the practical performance of the proposed framework under various scenarios. We generate 300 random instances, where $n = 30$ is the number of items, of which we select $q = 15$ items. We consider $K = 50$ scenarios and running the experiment for the  $p$-norm with $p\in\{1, 2, \ldots, 15\}$. The scenario cost values are chosen i.i.d. uniformly from ${1,\ldots,100}$.

\begin{figure}[htbp]
    \centering
    \includegraphics[height=7cm]
	{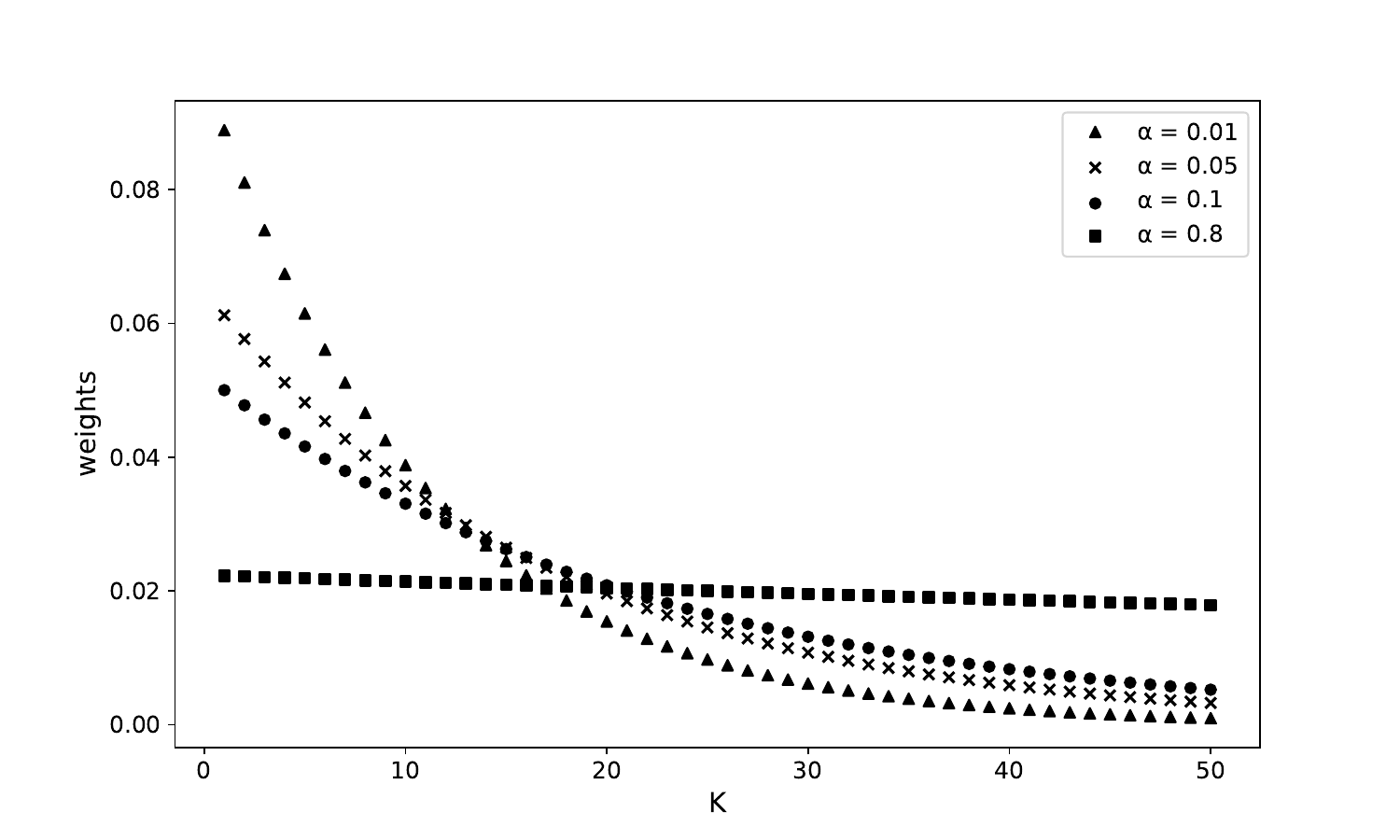}
    \caption{Preference vectors $\pmb{w}$ of dimension $K=50$ for different $\alpha$-values.
	}
    \label{fig:weight_generator}
\end{figure}

For the weight vector $\pmb{w}$, we use the generator functions
\[g_{\alpha}(z) = \frac{{1 - \alpha^z}}{{1 - \alpha}} \]
described in \cite{kasperski2016using} and define
\[w_k = g_{\alpha}\left(\frac{k}{K}\right) -g_{\alpha}\left(\frac{k-1}{K}\right) \quad k\in[K]. \]
The weight generator function allows us to control the distribution of weights in the OWA and OWAR criteria. By varying the value of $\alpha$, we can adjust the weights to be more conservative or more risk-neutral.$\alpha$ This flexibility enables us to explore a range of scenarios and assess the performance under different risk preferences. We consider two cases: one to model more risk-averse decision makers and another to model more risk-neutral decision makers, i.e., we use $\alpha=0.05$ (risk-averse) and $\alpha=0.8$ (risk-neutral) as our basic settings. The generated preference vectors are visualized in Figure~\ref{fig:weight_generator}.

\subsubsection{Results}\label{sec:results1}

\begin{figure}[H]
    \centering
    \begin{subfigure}[b]{0.495\textwidth}
        \includegraphics[width=\textwidth]{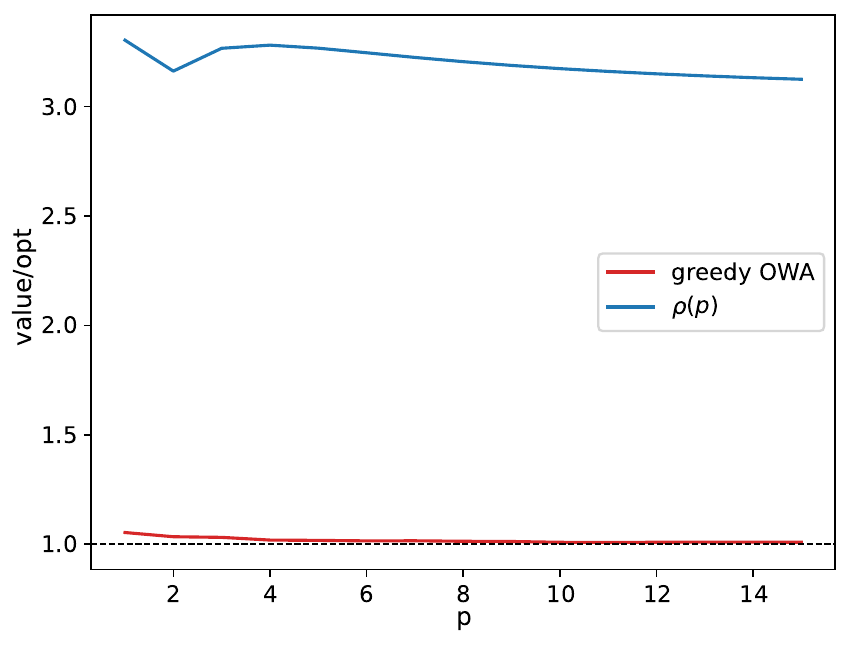}
        \caption{ $\alpha$=0.05}
    \end{subfigure}
    \hfill
    \begin{subfigure}[b]{0.495\textwidth}
        \includegraphics[width=\textwidth]{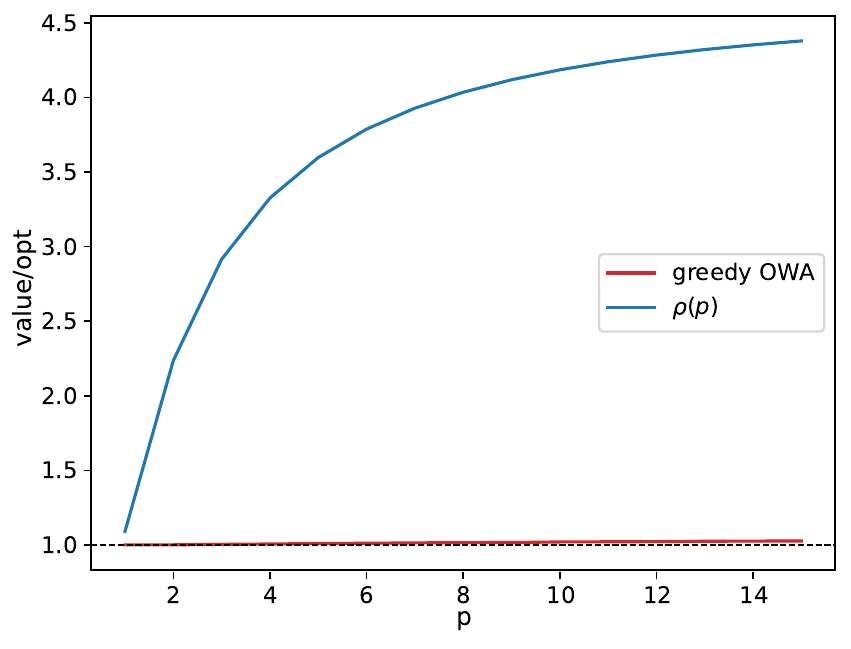}

        \caption{$\alpha$=0.8}
    \end{subfigure}
	\caption{Performance of greedy Algorithm for OWA.}\label{fig:GreedyOWA}
\end{figure}

\begin{figure}[H]
    \centering
    \begin{subfigure}[b]{0.49\textwidth}
        \includegraphics[width=\textwidth]{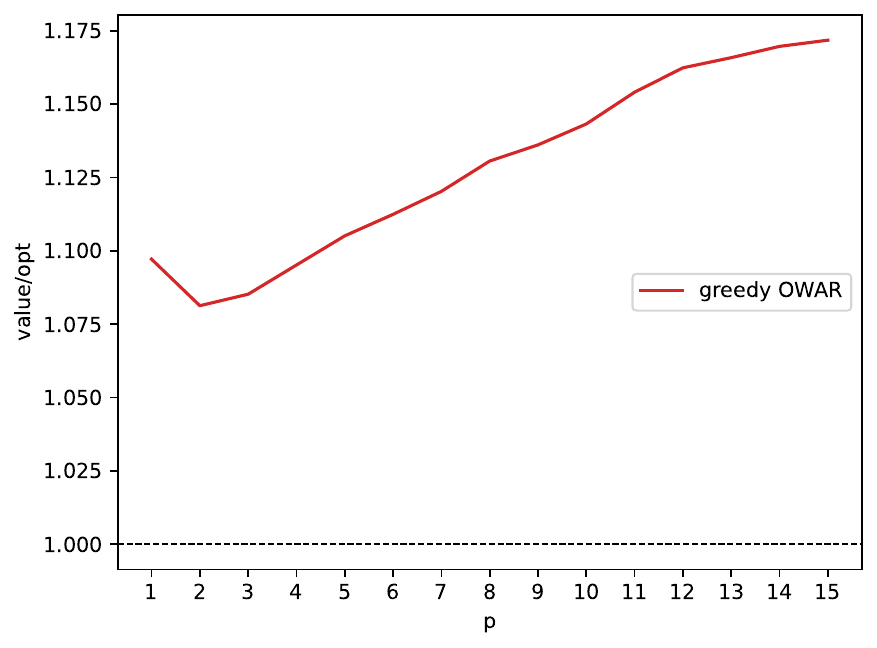}
        \caption{ $\alpha$=0.05}
    \end{subfigure}
    \hfill
    \begin{subfigure}[b]{0.49\textwidth}
        \includegraphics[width=\textwidth]{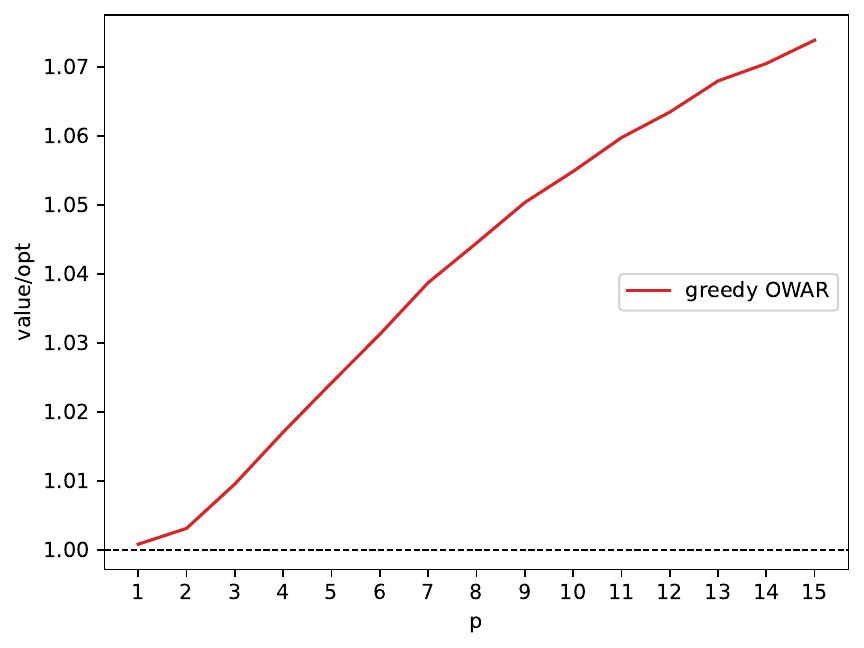}

        \caption{$\alpha$=0.8}
    \end{subfigure}
	\caption{Performance of greedy Algorithm for OWAR.}\label{fig:GreedyOWAR}
\end{figure}

In evaluating the performance of the greedy algorithm under the OWA criterion with $\pmb{b}=\pmb{0}$, we present the outcomes in Figure~\ref{fig:GreedyOWA}. Additionally, we showcase the results for the OWAR criterion with $\pmb{b}$, consisting of the optimal objective values over all scenarios in Figure \ref{fig:GreedyOWAR}. The greedy algorithm demonstrates commendable performance, exhibiting effectiveness, particularly with risk-neutral weights for lower $p$-norms. Moreover, it excels for higher $p$ values when confronted with risk-averse weights. The function $\rho(p)$, depicted similarly to Figure~\ref{fig2}, aligns with the theoretical expectations. Notably, as the $\alpha$ values approach 1, $\rho(p)$ tends to optimize for $p=1$.

Comparatively, the performance of the greedy algorithm in the context of OWAR (Figure~\ref{fig:GreedyOWAR}) exhibits a decline as $p$ increases. This trend arises from the algorithm's pursuit of regret minimization. Notably, the approximation guarantee for OWAR, which is noticeably less favorable than for OWA (refer to Theorem \ref{thmowa} and Theorem \ref{thmapprb}), is omitted from the plots due to its inferior performance. It is evident that the greedy solutions evaluated under OWA and OWAR objectives
 significantly outperform the theoretical guarantees across various setups, underscoring their efficacy in practical decision-making scenarios.

\subsection{Experiment 2: Performance of OWAR decisions}
\subsubsection{Setup}\label{sec:setup2}

In the second experiment, we investigate the performance of the ordered weighted averaging regret (OWAR) criterion in comparison to using the min-max regret or classic ordered weighted averaging (OWA) approach. We consider instances from the type of randomly generated selection problems.

To generate selection problems with $\X = \{ \pmb{x}\in\{0,1\}^n : \sum_{i\in[n]} x_i = q\}$, we focus on $n = 40$ items and $q = 20$ items to be selected. Each instance has $K=50$ scenarios, with each scenario value is chosen i.i.d. uniformly from $\{1,\ldots,100\}$. We generate 100 random instances this way.

Both OWA and OWAR require a weight vector $\pmb{w}$. To study the range from risk-averse to risk-neutral decision making we define vectors $\pmb{w}^k$ as $w^k_1=w^k_2=\ldots = w^k_k = 1/k$ and $w^k_{k+1}=w^k_{k+2} = \ldots = w^k_{K} = 0$, i.e., the first $k$ vector entries share the total weight uniformly. We generate vectors $\pmb{w}^k$ for $k=5,10,\ldots,50$, resulting in a total of 10 different weight vectors used to calculate the corresponding OWA and OWAR solutions.  We denote the resulting solutions as $\OWA_k$ and $\OWAR_k$ in the subsequent analysis.

For each instance we consider, we calculate a min-max regret solution and ten solutions $\OWA_k$ and $\OWAR_k$, respectively. To evaluate the quality of these solutions, we calculate the objective value of each solution in each of the 21 different decision criteria, resulting in a performance matrix.

\subsubsection{Results}

\begin{figure}[H]
	\begin{center}
	\includegraphics[width=1\textwidth]{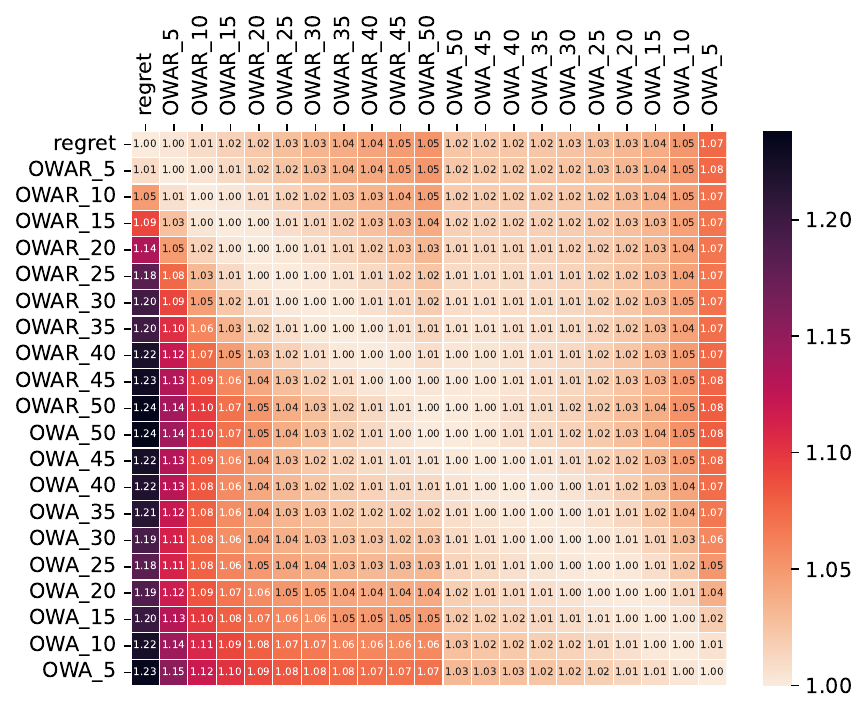}
	\end{center}
	\caption{Performance of OWAR decisions.}\label{fig:SeP}
	\end{figure}

We present a heat map in Figure~\ref{fig:SeP} to visualize the results in the case of the selection problem instances. The rows correspond to the decision criterion used to calculate a solution, while a column represents decision criterion used to calculate its objective value. All values are first normalized with respect to the best value of the column and then averaged over the 100 instances. By construction, there is a value of $1.00$ along the diagonal. As an example, the value $1.05$ in row ''regret'' and column ''$\OWAR_{45}$'' means that the min-max regret solution has an objective value with respect to $\OWAR_{45}$ that is on average $5\%$ higher than the optimal objective value for $\OWAR_{45}$. Note that $\OWAR_{50}$ and $\OWA_{50}$ yield the same optimal solution, rendering the two rows are identical (though the corresponding columns differ in objective values).

The heatmap illustrates that the OWAR criterion performs similarly to the min-max regret criterion for very conservative weights ($k=5$) and is equivalent to the OWA criterion for risk-neutral weights ($k=50$). In other words, the OWAR criterion provides a way to interpolate between the two extremes. Moreover, we have computed the average values of each row to assess the overall performance of each criterion (see Appendix~\ref{App::average}).
Upon analyzing the results, it is evident that the OWAR criterion consistently outperforms other criteria across different weight configurations. It consistently exhibits a lower average value, indicating superior performance in terms of minimizing regret. This finding aligns with the observation made in the heatmap, where OWAR is shown to be equivalent to the min-max regret criterion for conservative weights and equivalent to the OWA criterion for risk-neutral weights.

\subsection{Experiment 3: Performance of aggregation methods}
\subsubsection{Setup}\label{sec:setup3}

In the third experiment, we investigate the performance of scenario aggregation methods, specifically the $\ell$-aggregation and $K$-means algorithm, but also the guarantee given in Section~\ref{sec:aggregate}. We aim to assess the effectiveness of these algorithms in aggregating scenarios and understand their implications for decision-making under the OWAR criterion.

Similar to Experiment 1 and Experiment 2, we consider randomly generated selection problems with $\X = \{ \pmb{x}\in\{0,1\}^n : \sum_{i\in[n]} x_i = q\}$, where $n = 30$ items and $q = 15$ items to be selected. Each instance has $K=50$ scenarios, and the scenario value is chosen i.i.d. uniformly from ${1,\ldots,100}$. We generate 50 random instances using this method.

For the $\ell$-aggregation method, we vary the number of aggregated scenarios using $\ell\in\{1,2,5,10,25,50\}$.  We aggregate the weights with $\overline{w}_i = w_{(i-1)\ell+1} + \ldots + w_{i\ell}$, the scenario values (costs) using  $\overline{\pmb{c}}_i= \sum_{j\in[\ell]} \pmb{c}_{(i-1)\ell  + j}/\ell$ and the $\pmb{b}$ values using $\overline{b}_i=\sum_{j\in[\ell]} b_{(i-1)\ell  + j}/\ell$, as referred in Section~\ref{sec:aggregate}.
Moreover, we calculate $\phi$ for each aggregated solution, with $\phi = \max_{k\in[K/\ell]} \left(\sum_{i\in[k]} w_i / \sum_{i\in[k]} \overline{w}_i\right)$ and the guarantee with $\ell\phi$.

The $K$-means algorithm utilizes the well-known $K$-means clustering technique~(see, e.g.,~\cite{LRU19})
 to group scenarios and calculate aggregated weights.
 For the algorithm, we use the same range for the number of clusters $\{1,\ldots,50\}$. When aggregating the weights, we have to consider the case where the total number of scenarios is not divisible by the desired number of clusters. In this case we handle the remainder scenarios as follows: We start from the last cluster and move backwards, assigning one additional scenario to each cluster until no remainder is left. This ensures that all scenarios are included in the aggregation, and the additional scenarios are distributed as evenly as possible among the last clusters. This method of handling the remainder of the scenarios helps to maintain balance in the aggregation process while ensuring that all scenarios are taken into account.

We generate plots to provide visual representations of the experimental results obtained from the $\ell$-aggregation algorithm, the given guarantee, and the $K$-means algorithm for the selection problem when varying the $\alpha$ parameter.
We construct the plots such that the y-axis represents the ratio of OWAR achieved by each approach relative to the optimal solution. It quantifies the effectiveness of the algorithms in minimizing the overall OWAR, with a lower ratio indicating better performance. The $\ell$-aggregation is plotted on the x-axis for $k/\ell$, where $\ell$ is the number of aggregated scenarios and $k$ is the total number of scenarios.

\subsubsection{Results}\label{sec:results3}

\begin{figure}[H]
    \centering
    \begin{subfigure}[b]{0.495\textwidth}
        \includegraphics[width=\textwidth]{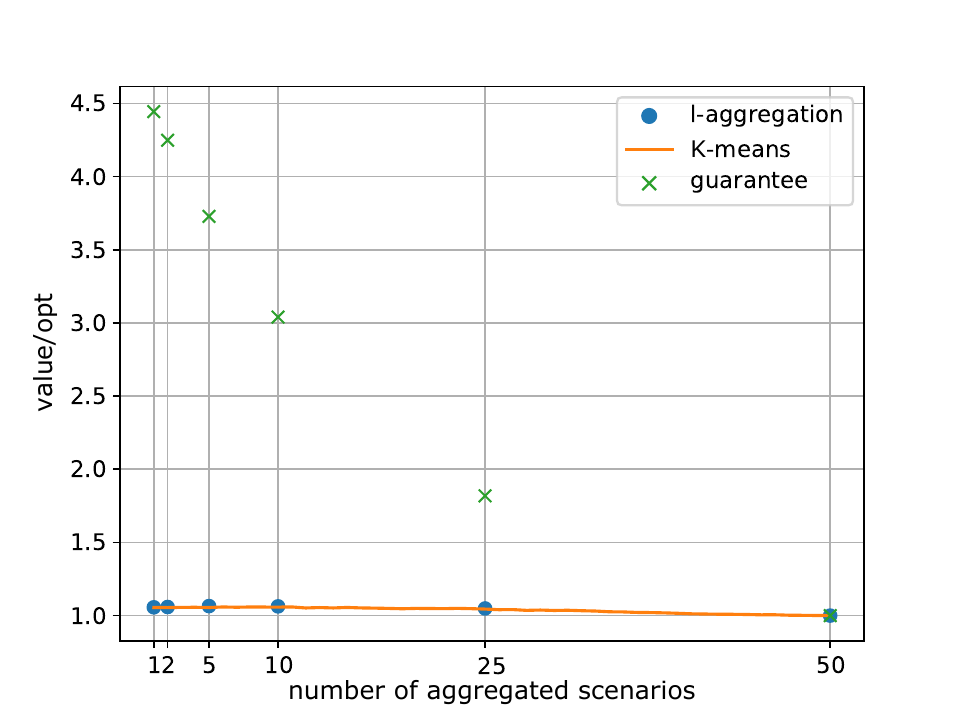}
        \caption{ $\alpha$=0.01}\label{fig:Agga}
    \end{subfigure}
    \hfill
    \begin{subfigure}[b]{0.495\textwidth}
        \includegraphics[width=\textwidth]{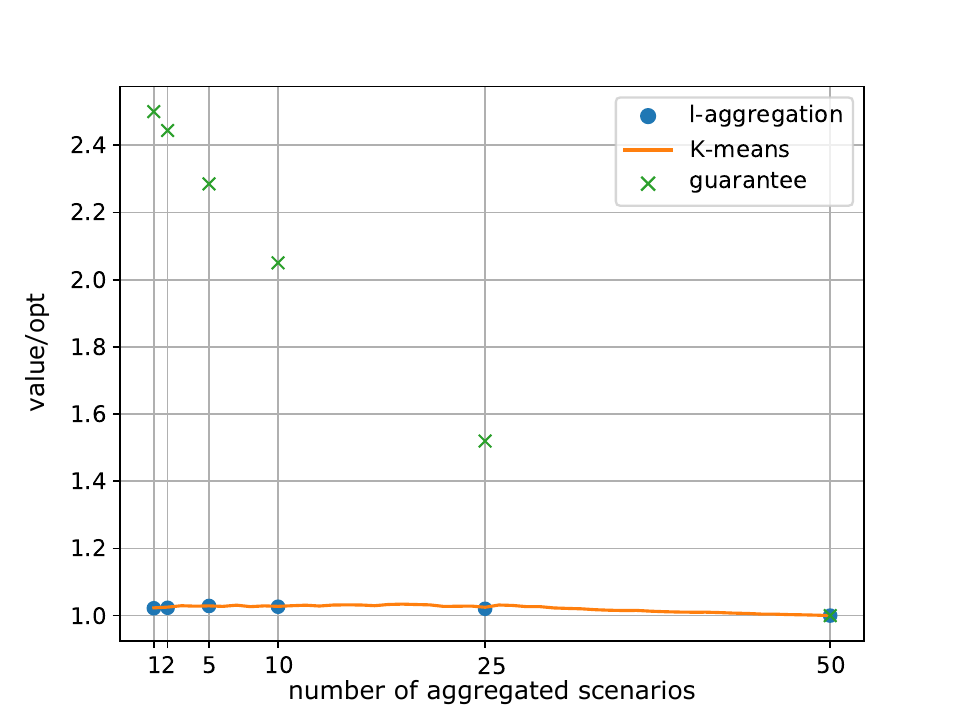}
        \caption{$\alpha$=0.1}\label{fig:Aggb}
    \end{subfigure}
	\caption{Performance of Aggregation Methods.}\label{fig:Agg}
\end{figure}

The results of the aggregation methods for the selection problem are presented in Figure~\ref{fig:Agg}. Both aggregation methods exhibit good performance, with the $K$-means algorithm slightly outperforming the $\ell$-aggregation method. Specifically, the $K$-means algorithm achieves a lower OWAR ratio than the $\ell$-aggregation method for all values of $\alpha=0.01$. On the other hand, the $\ell$-aggregation method achieves a lower OWAR ratio than the guarantee for all values of $\alpha$, except for $\ell=1$ and the case where the number of aggregated scenarios equals $K$. In this specific case, the ratio is 1 for all methods.

As expected, the conservative guarantee is outperformed by the $\ell$-aggregation method, indicating that the bounds could be tightened. In contrast to the results of a risk-averse aggregation (\ref{fig:Agga}), a clear trend is not observed for the more risk-neutral aggregation (\ref{fig:Aggb}). This lack of trend is due to the fact that risk-neutral aggregation is less sensitive to the number of aggregated scenarios compared to risk-averse aggregation. However, it is evident that the guarantees given in this case are bound tighter than in the risk-averse scenario.

\section{Conclusions}
\label{sec:conclusion}

In this paper, we have studied a class of optimization problems with uncertain objective functions. This uncertainty has been modeled using a discrete scenario set containing a finite number of cost scenarios. We have introduced an additional vector $\pmb{b}$ to modify the objective values under scenarios, and we have used OWA to aggregate the resulting vector of affine functions depending on scenarios and $\pmb{b}$. This approach allowed us to generalize both the robust min-max and min-max regret optimization.
For this general setting, several new complexity results have been provided. In particular, in the case of combinatorial optimization, we used norm-based estimates to find new approximation guarantees that improved the previously best-known results for the more specific case of classic OWA optimization.  We have shown general results demonstrating relationships between the OWA and the norm or reference point optimizations. We have applied them to particular problems for which some results in this area have been previously established.
In three computational experiments, we explored the practical implications of our theoretical developments. Experiment 1 spotlighted the effectiveness of the greedy algorithm for optimizing OWA and OWAR criteria.
 In Experiment 2, we investigated the OWAR criterion's performance relative to min-max regret and classic OWA, revealing its capacity to interpolate between conservative and risk-neutral decision-making. Experiment 3 turned attention to scenario aggregation methods, demonstrating that both $\ell$-aggregation and the $K$-means algorithm effectively minimize OWAR, with the latter exhibiting a slight edge.

There are still some open problems regarding the approximability of the OWA optimization for the class of combinatorial problems. Namely, there is still a gap between the positive and negative results shown in this paper for some classes of weight distributions - in particular, when $w_1=\Theta(\frac{1}{\sqrt{K}})$. Also, better approximation algorithms can be constructed for particular optimization problems by taking into account the inner structure of $\mathcal{X}$. In further research, the impact of the additional vector $\pmb{b}$ on the set of decision-maker preferences that can be modeled can also be investigated.

\section*{Acknowledgements} Marc Goerigk and Werner Baak were supported by the Deutsche Forschungsgemeinschaft (DFG) through grant 448792059. Adam Kasperski and Pawe{\l} Zieli\'nski were supported by the National Science Centre, Poland, grant 2022/45/B/HS4/00355.

\newpage
\appendix
\section{Averages of weight distributions}\label{App::average}

In Table~\ref{tab:combined_table}, we show the average performance (i.e., the average of each row) for the decision criteria presented in Figure~\ref{fig:SeP} for the selection problem.

\begin{table}[htbp]
	\centering
	\caption{Averaged performance for decision criteria.}
	\label{tab:combined_table}
	\begin{tabular}{c|c}
	  Criterion & Average performance \\
	  \hline
		regret & 1.030  \\
		$\OWAR_{5}$ & 1.029  \\
		$\OWAR_{10}$ & 1.027  \\
		$\OWAR_{15}$ & 1.026  \\
		$\OWAR_{20}$ & 1.027  \\
		$\OWAR_{25}$ & 1.029  \\
		$\OWAR_{30}$ & 1.031  \\
		$\OWAR_{35}$ & 1.032  \\
		$\OWAR_{40}$ & 1.036  \\
		$\OWAR_{45}$ & 1.040  \\
		$\OWAR_{50}$ & 1.045 \\
		$\OWA_{50}$ & 1.045  \\
		$\OWA_{45}$ & 1.040  \\
		$\OWA_{40}$ & 1.039  \\
		$\OWA_{35}$ & 1.038  \\
		$\OWA_{30}$ & 1.037  \\
		$\OWA_{25}$ & 1.039  \\
		$\OWA_{20}$ & 1.042  \\
		$\OWA_{15}$ & 1.048  \\
		$\OWA_{10}$ & 1.055  \\
		$\OWA_{5}$ & 1.063 \\
\end{tabular}
  \end{table}

\end{document}